\documentclass[11pt]{article}
\usepackage{latexsym,eucal,amsmath,amssymb,color}
\usepackage{pstricks}
\usepackage{setspace}

\newrgbcolor{VertexColour}{1 1 1}
\newrgbcolor{EdgeColour}{0 0 0}

\newtheorem{Theorem}{Theorem}[section]
\newtheorem{Corollary}{Corollary}[section]

\newenvironment{Proof}{\begin{trivlist} \item[] {\bf Proof.}}{\hfill $\Box$\end{trivlist}}

\newcommand{\sB}{\EuScript{B}}

\newcommand{\myK}{\mathcal K}

\renewcommand{\geq}{\geqslant}
\renewcommand{\leq}{\leqslant}

\title{Colourings of star systems}

\author{Iren Darijani and  David A. Pike\\
Department of Mathematics and Statistics\\
Memorial University of Newfoundland\\
St.\ John's, NL, Canada. A1C 5S7}

\begin{document}
	
\maketitle

\abstract 
An $e$-star is a complete bipartite graph $K_{1,e}$. An $e$-star system of order $n>1$, $S_e(n)$, is a partition of the edges of the complete graph $K_n$ into $e$-stars. An $e$-star system is said to be $k$-colourable if its vertex set can be partitioned into $k$ sets (called colour classes) such that no $e$-star is monochromatic. The system $S_e(n)$ is $k$-chromatic if $S_e(n)$ is $k$-colourable but is not $(k-1)$-colourable. If every $k$-colouring of an $e$-star system can be obtained from some $k$-colouring $\phi$ by a permutation of the colours, we say that the system is uniquely $k$-colourable. In this paper, we first show that for any integer $k\geq 2$, there exists a $k$-chromatic 3-star system of order $n$ for all sufficiently large admissible $n$. Next, we generalize this result for $e$-star systems for any $e\geq 3$. We show that for all $k\geq 2$ and $e\geq 3$, there exists a $k$-chromatic $e$-star system of order $n$ for all sufficiently large $n$ such that $n\equiv 0,1$ (mod $2e$). Finally, we prove that for all $k\geq 2$ and $e\geq 3$, there exists a uniquely $k$-chromatic $e$-star system of order $n$ for all sufficiently large $n$ such that  $n\equiv 0,1$ (mod $2e$).

\vspace*{\baselineskip}

\vspace*{\baselineskip}

\begin{doublespace}	

\section{Introduction}
A {\em $G$-decomposition} of a graph $H$ is a pair $(V,\mathcal{B})$ where $V$ is the set of vertices of $K$ and $\mathcal{B}$ is a set of subgraphs of $K$, each isomorphic to $G$, whose edge sets partition the edge set of $K$. A {\em $G$-design} of order $n$ is a $G$-decomposition of the complete graph $K_n$ on $n$ vertices. A {\em complete graph} is a simple graph in which every pair of distinct vertices is connected by a unique edge. A $G$-design $(V,\mathcal{B})$  is said to be {\em weakly $k$-colourable} if its vertex set can be partitioned into $k$ sets (called colour classes) such that no subgraph belonging to $\sB$ is monochromatic. The  $G$-design is {\em $k$-chromatic} if it is $k$-colourable but is not $(k-1)$-colourable. If a $G$-design is $k$-chromatic, we say that its chromatic number is $k$. A colouring of a  $G$-design is said to be {\em equitable} if the cardinalities of the colour classes differ by at most one. It is {\em strongly equitable} if the colour classes are of the same size. A $G$-design will be called {\em (strongly) equitably $k$-chromatic} if it is $k$-chromatic and admits a (strongly) equitable $k$-colouring. If every $k$-colouring of a $G$-design $(V,\mathcal B)$ can be obtained from some $k$-colouring $\phi$ by a permutation of the colours, we say that $(V,\mathcal B)$ has a {\em unique} $k$-colouring, or that $(V,\mathcal B)$ is {\em uniquely $k$-colourable}.

Weak colourings of $G$-designs were first studied for $K_3$-designs, i.e., Steiner triple systems. In particular, de Brandes, Phelps and R\"{o}dl \cite{de Brandes} proved that for any $k\geq 3$, there exists some integer $v_k$ such that for all admissible $v\geq v_k$, there is a $k$-chromatic Steiner triple system of order $v$. Moreover, in \cite{Burgess} Burgess and Pike showed that for all $k\geq 2$ and even $m\geq 4$ there exists a $k$-chromatic $m$-cycle system. An $m$-cycle system of order $n>1$ is a partition of the edges of the complete graph $K_n$ into $m$-cycles. Also, Horsley and Pike showed that for all $k\geq 2$ and $m\geq 3$ with $(k,m)\neq (2,3)$, there exist $k$-chromatic $m$-cycle systems of all admissible orders greater than or equal to some integer $n_{k,m}$ \cite{Horsley}. 

A bipartite graph $G_{m,n}$ is a graph whose vertex set can be partitioned into two subsets $V_1$ and $V_2$ with $m$ and $n$ vertices each, such that every edge of $G_{m,n}$ has one vertex in $V_1$ and one vertex in $V_2$. If $G_{m,n}$ contains every edge joining $V_1$ and $V_2$, then it is called a {\em complete bipartite graph} and denoted by $K_{m,n}$. For an integer $e\geq 1$, an {\em $e$-star} is a complete bipartite graph $K_{1,e}$.  An $e$-star with star centre $x$ and edges $\{x,a_1\}$, $\{x,a_2\}$, $\ldots$, $\{x,a_e\}$ is denoted by either  $\{x;a_1,\ldots,a_e\}$ or $\{x;A\}$ where $A=\{a_1,\ldots,a_e\}$. When $G$ is an $e$-star, a $G$-design of order $n$ is called an {\em $e$-star system} of order $n$ and denoted by $S_e(n)$. The necessary and sufficient conditions for the existence of an $e$-star system of order $n$ are that $e$ divides $\frac{n(n-1)}{2}$ and $n \geq 2e$~\cite{Yamamoto1975}. A positive integer $n$ is said to be {\em admissible} if there exists an $e$-star system of order $n$. 

A {\em $G$-factor} of a graph $H$ is a spanning subgraph of $H$, each component of which is isomorphic to $G$. If the edges of $H$ can be partitioned into $G$-factors, then we say $H$ has a {\em $G$-factorization}. The existence of a $K_{1,e}$-factorization of $K_n$ was studied and completely settled in 1993~\cite{Yu}. If $e$ is even, a $K_{1,e}$-factorization of $K_n$ exists if and only if $n \equiv 1$ (mod $2e$) and $n \equiv 0$ (mod $e+1$). If $e$ is odd, there does not exist a $K_{1,e}$-factorization of $K_n$ for any $n$.

In this paper, we investigate $k$-colourings of $e$-star systems for $e\geq 3$. Note that a $3$-star system of order $n$ exists if and only if $n\equiv 0,1$ (mod 3) and $n\geq 6$. In Section~2, we first construct equitably 2-chromatic 3-star systems of all admissible orders. We then show that if there exists a $k$-chromatic $3$-star system of order $n_0$, then there exists a $k$-chromatic $3$-star system of order $n$ for all admissible $n>n_0$ and $k\geq 2$. Next, from a $(k-1)$-chromatic 3-star system, we construct a $k$-chromatic 3-star system for all $k \geq 3$. Finally, we finish this section by showing that for any integer $k \geq 2$, there exists some integer $n_k$ such that for all admissible $n \geq n_k$, there exists a $k$-chromatic $3$-star system of order $n$.

 In Section~3, we generalize the results in Section~2 for $e$-star systems for all $e\geq 3$. We first construct equitably 2-chromatic $e$-star systems of order $2e$ for all $e\geq 3$. We then show that if there exists a $k$-chromatic $e$-star system of order $n_0$ such that $n_0\equiv 0,1$ (mod $2e$), then there exists a $k$-chromatic $e$-star system for all $n >n_0$ such that $n\equiv 0,1$ (mod $2e$) for all $k\geq 2$ and $e\geq 4$. Next, from a $(k-1)$-chromatic $e$-star system of order $n_{k-1}\equiv 0$ (mod $2e$), we construct a $k$-chromatic $e$-star system of order $n_k \equiv 0$ (mod $2e$) for all $k \geq 3$. Finally, we finish this section by showing that for any integer $k \geq 2$, there exists some integer $n_k$ where $n_k \equiv 0$ (mod $2e$) such that for all $n \geq n_k$ where $n\equiv 0,1$ (mod $2e$), there exists a $k$-chromatic $e$-star system of order $n$.

We also study unique colourings for $e$-star systems. The idea of uniquely colourable designs has already arisen in the content of Steiner triple systems. In 2003, Forbes showed that for every admissible $v\geq 25$, there exists a 3-balanced Steiner triple system with a unique 3-colouring and also a Steiner triple system  which has a unique, nonequitable 3-colouring~\cite{Forbes}, where a Steiner triple system is said to be {\em 3-balanced} if every 3-colouring of it is equitable. In Section~4, we first construct a strongly equitable uniquely 2-chromatic $e$-star system. Next, we construct a strongly equitable $k$-chromatic $e$-star system from a strongly equitable uniquely $(k-1)$-chromatic $e$-star system. We then show how to construct a strongly equitable uniquely $k$-chromatic $e$-star system from a strongly equitable $k$-chromatic $e$-star system. Finally, we finish this section by showing that for any integer $k\geq 2$, there exists some integer $n_k$ where $n_k\equiv 0$ (mod $2e$) such that for all $n \geq n_k$ where $n\equiv 0,1$ (mod $2e$), there exists a uniquely $k$-chromatic $e$-star system of order $n$.

\section{$\mathbf k$-colourings of 3-star systems}
We initially concentrate on 3-star systems and we show that for any integer $k \geq 2$, there exists some integer $n_k$ such that for all admissible $n \geq n_k$, there exists a $k$-chromatic $3$-star system of order $n$. We first construct an equitably 2-chromatic 3-star system of order $n$ for all admissible $n$.

\begin{Theorem}\label{thm2}
For each admissible order $n$, there exists an equitably 2-chromatic 3-star system of order $n$.
\end{Theorem}

\begin{Proof}
We first construct an equitably 2-chromatic 3-star system $S_3(6)$. Let $V=\{1,2,3,4,5,6\}$, $\mathcal B=\big \{ \{1;3,5,6\}$, $\{2;1,3,6\}$, $\{4;1,2,3\}$, $\{5;2,3,4\}$, $\{6;3,4,5\}\big\}$, $R=\{1,3,5\}$, and $Y=\{2,4,6\}$. Then $(V,\mathcal B)$ is an equitably 2-chromatic 3-star system of order six with colour classes $R$ and $Y$.

Now, suppose that there exists an equitably  2-chromatic 3-star system $S_3(3t)$, $(V,\mathcal B)$, where $t\geq 2$, $V=\{1,\ldots,3t\}$ is the set of points which is partitioned into two subsets $R$ and $Y$ and $\mathcal B$ is the set of blocks. Without loss of generality, if $3t$ is odd, we can let $R=\{1,3,5,\ldots,3t\}$, $Y=\{2,4,6,\ldots,3t-1\}$ and if $3t$ is even, $R=\{1,3,5,\ldots,3t-1\}$, $Y=\{2,4,6,\ldots,3t\}$. 

Next, we construct an equitably 2-chromatic 3-star system $S_3(3t+1)$, $(\hat{V},\hat{\mathcal B})$, from $(V,\mathcal B)$. Let $\hat{V}=V\cup\{3t+1\}$ and $\hat{\mathcal B}=\mathcal B \cup \big\{\{3t+1;1,2,3\},\{3t+1;4,5,6\},\ldots,\{3t+1;3t-2,3t-1,3t\}\big\}$. Let $R=\{1,3,5,\ldots,3t\}$, $Y=\{2,4,6,\ldots,3t+1\}$ if $3t$ is odd and $R=\{1,3,5,\ldots,3t+1\}$, $Y=\{2,4,6,\ldots,3t\}$ if $3t$ is even. Then $(\hat{V},\hat{\mathcal B})$ is an equitably 2-chromatic 3-star system of order $3t+1$ with colour classes $R$ and $Y$. 

Finally, we construct an equitably 2-chromatic 3-star system $S_3(3t+3)$, $(\tilde{V},\tilde{\mathcal B})$, from $(V,\mathcal B)$. Let $\tilde{V}=V\cup\{3t+1,3t+2,3t+3\}$ and $\tilde{\mathcal B}=\mathcal B\cup\big\{\{3t+1;1,2,3\},\ldots,\{3t+1;3t-5,3t-4,3t-3\},\{3t+2;1,2,3\},\ldots,\{3t+2;3t-5,3t-4,3t-3\},\{3t+3;1,2,3\},\ldots,\{3t+3;3t-5,3t-4,3t-3\},\{3t+1;3t-2,3t-1,3t+2\},\{3t+2;3t-2,3t-1,3t+3\},\{3t+3;3t-2,3t-1,3t+1\},\{3t;3t+1,3t+2,3t+3\}\big\}$. Let $R$ be the set of odd elements of $\tilde{V}$ and $Y$ be the set of even elements of $\tilde{V}$. Then $(\tilde{V},\tilde{\mathcal B})$ is an equitably 2-chromatic 3-star system of order $3t+3$ with colour classes $R$ and $Y$.  
\end{Proof}
We now show how to construct a $k$-chromatic 3-star system from a smaller $k$-chromatic 3-star system.
\begin{Theorem}\label{thm1}
Let $k\geq 2$. If there exists a $k$-chromatic $3$-star system of order $n_0$, then there exists a $k$-chromatic $3$-star system of order $n$ for all admissible $n>n_0$. 
\end{Theorem}

\begin{Proof}
Suppose that there exists a $k$-chromatic $3$-star system of order $n_0$, $(V,\mathcal B)$, where $V=\{1,\ldots, n_0\}$ is the set of points and $\mathcal B$ is the set of blocks. 
 Given a $k$-colouring of $(V,\mathcal B)$ with colours $1,2,\ldots,k$ and colour classes $C_1,C_2,\ldots,C_k$, let $R=\bigcup\limits_{i=1}^{\ell}C_i$ and $Y=\bigcup\limits_{i=\ell+1}^{k}C_i$ for some integer $\ell$ such that $1\leq \ell < k$. 
Let $r_1, r_2,\ldots,r_{|R|}$ be the elements of $R$ and $y_1, y_2, \ldots, y_{|Y|}$ be the elements of $Y$. Observe that $R$ and $Y$ form a partition of $V$. Without loss of generality assume that $|R|\geq |Y|$. 

\textbf {Case 1.} Suppose that $n_0\equiv 0$ (mod 3). Then $n_0=3t$, $t\geq 2$.

 First, we construct a $k$-chromatic $S_3(3t+1)$, $(\hat{V}, \hat{\mathcal B})$, from $(V,\mathcal B)$. Let $\hat{V}=V\cup \{3t+1\}$.

If $|R|=|Y|$, let
$\hat{\mathcal B}=\mathcal B \cup \mathcal T$ where $\mathcal T$ is the set
$\big \{\{3t+1;r_1,r_2,y_1\}$, $\{3t+1;y_2,y_3,r_3\}$, $\ldots,$ $\{3t+1;r_{|R|-2},r_{|R|-1},y_{|Y|-2}\},$ $\{3t+1;y_{|Y|-1},y_{|Y|},r_{|R|}\}\big \}.$ 

Otherwise, $|R|>|Y|$. If $|Y|$ is even then let $m=\frac{|Y|}{2}$ and $\hat{\mathcal B}=\mathcal B \cup \mathcal T$ where $\mathcal T$ is the set
$\big\{\{3t+1;y_1,y_2,r_1\}$, $\ldots,$ $\{3t+1;y_{|Y|-1},y_{|Y|},r_m\}$, $\{3t+1;r_{m+1},r_{m+2},r_{m+3}\}$, $\ldots,$ $\{3t+1;r_{|R|-2},r_{|R|-1},r_{|R|}\}\big\}$.
If $|Y|$ is odd then let $m=\frac {|Y|+3}{2}$ and 
$\hat{\mathcal B}=\mathcal B \cup \mathcal T$ where $\mathcal T$ is the set 
$\big\{\{3t+1;y_1,y_2,r_1\}$, $\ldots,$ $\{3t+1;y_{|Y|-2},y_{|Y|-1},r_{m-2}\},$ $\{3t+1;y_{|Y|},r_{m-1},r_m\},$ $\{3t+1;r_{m+1},r_{m+2},r_{m+3}\},$ $\ldots,$
$\{3t+1;r_{|R|-2},r_{|R|-1},r_{|R|}\}\big\}$. 

Note that $(\hat{V}, \hat{\mathcal B})$ is not $(k-1)$-colourable because it contains a copy of $(V,\mathcal B)$. Observe that $C_1,\ldots,C_{k-1},C_k\cup \{3t+1\}$ are the colour classes of a $k$-colouring of $(\hat{V},\hat {\mathcal B})$. Therefore, $(\hat{V},\hat {\mathcal B})$ is a $k$-chromatic $3$-star system of order $3t+1$.

Next, we construct a $k$-chromatic $S_3(3t+3)$, $(\tilde{V},\tilde{\mathcal B})$, from $(V,\mathcal B)$. Let $\tilde{V}=V\cup \{3t+1,3t+2,3t+3\}$ and $R'=R \setminus \{r_{|R|-2},r_{|R|-1},r_{|R|}\}$. Note that $|R| \geq 3$ since $n_0\geq 6$ and $|R| \geq |Y|$.

Let $\mathcal T$ be the set $\big \{\{3t+1;r_{|R|-2},r_{|R|-1},3t+2\},$ $\{3t+2;r_{|R|-2},r_{|R|-1},3t+3\},$ $\{3t+3,r_{|R|-2},r_{|R|-1},$ $3t+1\}$, $\{r_{|R|};3t+1,3t+2,3t+3\}\big \}$.

If $|R'|=|Y|$, let $\tilde {\mathcal B}=\mathcal B \cup \mathcal T \cup (\bigcup \limits_{i=1}^{i=3} \mathcal T_i)$ where $\mathcal T_i$ is the set $\big\{\{3t+i;r_1,r_2,y_1\},$ $\{3t+i;y_2,y_3,r_3\},$ $\ldots,$ $\{3t+i;r_{|R|-5},r_{|R|-4},y_{|Y|-2}\},$ $\{3t+i;y_{|Y|-1},y_{|Y|},r_{|R|-3}\}\big\}$ for $i \in \{1,2,3\}$.


If $|R'|>|Y|$ and $|Y|$ is even, then let $m=\frac{|Y|}{2}$ and $\tilde {\mathcal B}=\mathcal B \cup \mathcal T \cup (\bigcup \limits_{i=1}^{i=3} \mathcal T_i)$ where $\mathcal T_i$ is the set
$\big \{\{3t+i;y_1,y_2,r_1\},$ $\ldots,$ $\{3t+i;y_{|Y|-1},y_{|Y|},r_m\},$ $\{3t+i;r_{m+1},r_{m+2},r_{m+3}\},$ $\ldots,$ $\{3t+i;r_{|R|-5},r_{|R|-4},$ $r_{|R|-3}\}\big \}$  for $i \in \{1,2,3\}$.
 The edges can be decomposed into 3-stars in a similar manner when $|R'|>|Y|$ and $|Y|$ is odd.

Otherwise, $|R'|<|Y|$. Then either $|R'|=|Y|-1$ or $|R'|=|Y|-2$ or $|R'|=|Y|-3$.
If $|R'|=0$, then $|R|=|Y|=3$. Without loss of generalty, we may assume that $Y=\{1,2,3\}$ and $R=\{4,5,6\}$. Now let $\tilde {\mathcal B}=\mathcal B \cup \mathcal T_0$ where $\mathcal T_0$ is the set $\big\{ \{7;1,2,4\},$ $\{8;1,2,4\},$ $\{9;1,2,4\},$ $\{7;8,3,5\},$ $\{8;9,3,5\},$ $\{9;7,3,5\},$ $\{6;7,8,9\} \big\}$.
 If $|R'|\geq 1$ and $|R'|=|Y|-1$, let $\tilde {\mathcal B}=\mathcal B \cup \mathcal T \cup (\bigcup \limits_{i=1}^{i=3} \mathcal T_i)$ where $\mathcal T_i$ is the set
$\big\{\{3t+i;r_1,y_1,y_2\}, \big \{\{3t+i;r_2,r_3,y_3\},$ $\{3t+i;y_4,y_5,r_4\},$ $\ldots,$ $\{3t+i;r_{|R|-2},r_{|R|-1},y_{|Y|-2}\},$ $\{3t+i;y_{|Y|-1},y_{|Y|},r_{|R|}\}\big\},$ 
for $i \in \{1,2,3\}$.
 The edges can be decomposed into 3-stars in a similar manner when $1\leq |R'|<|Y|$ and $|R'|=|Y|-j$ for $j\in\{2,3\}$.
 
 Note that $(\tilde{V}, \tilde{\mathcal B})$ is not $(k-1)$-colourable because it contains a copy of $(V,\mathcal B)$. Observe that $C_1,\ldots,C_{k-1},C_k\cup \{3t+1,3t+2,3t+3\}$ are the colour classes of a $k$-colouring of $(\tilde{V},\tilde {\mathcal B})$. Therefore, $(\tilde{V},\tilde {\mathcal B})$ is a $k$-chromatic $3$-star system of order $3t+3$.
 
 \textbf {Case 2.} Suppose that $n_0 \equiv 1$ (mod 3). Then $n_0=3t+1$, $t\geq 2$. We construct a $k$-chromatic $S_3(3t+3)$, $(\hat{V},\hat {\mathcal B})$, from $(V,\mathcal B)$. Let $\hat{V}=V\cup \{3t+2,3t+3\}$. Note that $|R|\geq 4$ since $n_0\geq 7$ and $|R|\geq |Y|$. So without loss of generality we may assume that  $\{3t,3t+1\}\subset R$. Moreover, since the edge $\{3t,3t+1\}$ must be in some 3-star, we may also assume that there exists a 3-star $S=\{3t+1;3t,3t-1,3t-2\} \in \mathcal B$.  
 We will dismantle the $3$-star $S$ in each of the following subcases:
 
 \textbf {Case 2.1.} All the vertices of the $3$-star $S$ belong to the set $R$. Hence $R$ is the union of $\ell \geq 2$ distinct colour classes. Let $R'=R\setminus \{3t+1,3t,3t-1,3t-2\}$.

 Let $\mathcal T$ be the set $\big\{\{3t+2;3t+1,3t,3t-1\},$ $\{3t+3;3t+2,3t,3t-1\},$ $\{3t-2;3t+2,3t+1,3t+3\}$, $\{3t+1;3t+3,3t,3t-1\}\big\}$.
 
 If $|R'|=|Y|$, let $\hat {\mathcal B}=(\mathcal B\setminus\{S\}) \cup \mathcal T \cup (\bigcup \limits_{i=1}^{i=2} \mathcal T_i)$ where $\mathcal T_i$ is the set 
  $\big\{\{3t+i;r_1,r_2,y_1\},$ $\{3t+i;y_2,y_3,r_3\},$ $\ldots,$ $\{3t+i;r_{|R|-6},r_{|R|-5},y_{|Y|-2}\},$ $\{3t+i;y_{|Y|-1},y_{|Y|},r_{|R|-4}\}\big\}$ for $i\in \{1,2\}$. 
  
  If $|R'|>|Y|$ and $|Y|$ is even, then let $m=\frac {|Y|}{2}$ and $\hat {\mathcal B}=(\mathcal B\setminus\{S\}) \cup \mathcal T \cup (\bigcup \limits_{i=1}^{i=2} \mathcal T_i)$  where $\mathcal T_i$ is the set  $\big\{\{3t+i;y_1,y_2,r_1\},$ $\ldots,$ $\{3t+i;y_{|Y|-1},y_{|Y|},r_m\},$ $\{3t+i;r_{m+1},r_{m+2},r_{m+3}\},$ $\ldots,$ $\{3t+i;r_{|R|-6},r_{|R|-5},$ $r_{|R|-4}\}\big\}$
  for $i \in \{1,2\}$.
  The edges can be decomposed into 3-stars in a similar manner when  $|R'|>|Y|$ and $|Y|$ is odd.
  
  Otherwise $|R'|<|Y|$. Then either $|R'|=|Y|-1$ or $|R'|=|Y|-2$ or $|R'|=|Y|-3$ or $|R'|=|Y|-4$. If $|R'|=0$, then $Y=\{1,2,3\}$ and $R=\{4,5,6,7\}$. Let 
  $\hat {\mathcal B}=\mathcal (\mathcal B\setminus\{S\}) \cup \mathcal T_0$ where $\mathcal T_0$ is the set 
  $\big\{ \{8;1,2,3\}$, $\{9;1,2,3\}$, $\{8;7,6,5\}$, $\{9;8,6,5\}$, $\{4;8,7,9\}$, $\{7;9,6,5\} \big\}.$  If $|R'|\geq 1$ and  $|R'|=|Y|-4$, let $\hat {\mathcal B}=(\mathcal B\setminus\{S\}) \cup \mathcal T \cup (\bigcup \limits_{i=1}^{i=2} \mathcal T_i)$ where $\mathcal T_i$ is the set   $\big\{\{3t+i;r_1,y_1,y_2\},$ $\{3t+i;r_2,y_3,y_4\},$  $\{3t+i;r_3,y_5,y_6\},$ $\{3t+i;r_4,y_7,y_8\},$ $\{3t+i;r_5,r_6,y_9\},$ $\{3t+i;y_{10},y_{11},r_7\},$ $\ldots,$ $\{3t+i;r_{|R|-6},r_{|R|-5},y_{|Y|-2}\},$ $\{3t+i;y_{|Y|-1},y_{|Y|},r_{|R|-4}\}\big\}$      
  for $i \in \{1,2\}$.
 The edges can be decomposed into 3-stars in a similar manner when $1\leq |R'|<|Y|$ and $|R'|=|Y|-j$, $1\leq j\leq 3$.
  
  Note that $(\hat{V}, \hat{\mathcal B})$ is not $(k-1)$-colourable because it contains a copy of $(V,\mathcal B)$. Observe that $C_1,\ldots,C_{k-1},C_k\cup \{3t+2,3t+3\}$ are the colour classes of a $k$-colouring of $(\hat{V},\hat {\mathcal B})$. Therefore, $(\hat{V},\hat {\mathcal B})$ is a $k$-chromatic $3$-star system of order $3t+3$.
  
  \textbf {Case 2.2.} Exactly three vertices of the $3$-star $S$ belong to the set $R$. The edges can be decomposed into $3$-stars in a manner similar to the last case.

  \textbf {Case 2.3.} Exactly two vertices of the $3$-star $S$ belong to the set $R$. Let $
 R'=R \setminus \{3t,3t+1\}$ and $Y'=Y \setminus \{3t-1,3t-2\}$. 
   
  If $|R'|>|Y'|$ and $|Y'|$ is even, then let $m=\frac{|Y'|}{2}$ and  $\hat {\mathcal B}=(\mathcal B\setminus\{S\}) \cup \mathcal T \cup (\bigcup \limits_{i=1}^{i=2} \mathcal T_i)$  where $\mathcal T_i$ is the set 
  $\big\{\{3t+i;y_1,y_2,r_1\},$ $\ldots,$ $\{3t+i;y_{|Y'|-3},y_{|Y'|-2},r_m\},$ $\{3t+i;r_{m+1},r_{m+2},r_{m+3}\},$ $\ldots,$ $\{3t+i;r_{|R|-4},r_{|R|-3},$ $r_{|R|-2}\}\big\}$ for $i \in \{1,2\}$
 and $\mathcal T$ is the set
  $\big\{\{3t+2;3t+1,3t,3t-1\},$ $\{3t+3;3t+2,3t,3t-1\},$ $\{3t-2;3t+2,3t+1,3t+3\}$, $\{3t+1;3t+3,3t,3t-1\}\big\}$.
   The edges can be decomposed into 3-stars in a similar manner when  $|R'|>|Y'|$ and $|Y'|$ is odd; likewise when $|R'|=|Y'|$.

    Note that $(\hat{V}, \hat{\mathcal B})$ is not $(k-1)$-colourable because it contains a copy of $(V,\mathcal B)$. Observe that $C_1,\ldots,C_{k-1},C_k\cup \{3t+2,3t+3\}$ are the colour classes of a $k$-colouring of $(\hat{V},\hat {\mathcal B})$. Therefore, $(\hat{V},\hat {\mathcal B})$ is a $k$-chromatic $3$-star system of order $3t+3$. 
\end{Proof}

We now show how to iteratively construct a $k$-chromatic 3-star system from a $(k-1)$-chromatic 3-star system.
\begin{Theorem}\label{thm3}
Let $k\geq 3$. If there exists a $(k-1)$-chromatic 3-star system then there exists a $k$-chromatic 3-star system.
\end{Theorem}
\begin{Proof}
Given that there exists a $(k-1)$-chromatic 3-star system, then by Theorem~\ref{thm1}, there exists a $(k-1)$-chromatic 3-star system of order $n_{k-1}$ for some $n_{k-1} \in \mathbb Z$ such that $n_{k-1} \equiv 0$ (mod 3). Let $(U_0,\mathcal A_0)$ be a $(k-1)$-chromatic 3-star system of order $n_{k-1}$ with colour classes $C_1,\ldots,C_{k-1}$ such that $|C_1| \leq \cdots \leq |C_{k-1}|$ and each vertex of $C_s$ has colour $s$, $1\leq s\leq k-1$. For each $s\in \{1,\ldots,k-1\}$, let $C_s=\{c_1^s,\ldots,c_{|C_s|}^s\}$. For a positive integer $\ell$ that will be fixed later and for each $i\in \{1,\ldots,\ell\}$, let $U_i=U_0 \times \{i\}$ and $\mathcal A_i = \mathcal A_0\times \{i\}$ where $\mathcal A_0\times \{i\}$ denotes the set $\{S\times \{i\} \mid S\in \mathcal A_0\}$ and when $S=\{x;a,b,c\}$, $S\times \{i\}$ denotes $\big\{x\times \{i\};a\times \{i\},b\times \{i\},c\times \{i\} \big\}$. So $(U_i,\mathcal A_i)$ has colour classes $C_1\times\{i\},\ldots,C_{k-1}\times\{i\}$. If $2k \equiv 0$ (mod 3) then let $V=\{1,2,3,\ldots,2k\}$; otherwise let $V=\{1,2,3,\ldots,2k-1\}$.  Also let $U=\bigcup\limits_{i=1}^{\ell} U_i$ such that $V\cap U=\emptyset$. We will embed $(U_1,\mathcal A_1),\ldots,(U_{\ell},\mathcal A_{\ell})$ into a $k$-chromatic $3$-star system $(\hat V,\hat {\mathcal B})$ where $\hat V=V \cup U$. Since the edges of the complete graph on the set $V$ admit a decomposition into 3-stars, we let $(V,\mathcal B)$ be an arbitrary 3-star system of order $2k-1$. We now need to decompose the edges between $V$ and $U$ and the edges between $U_i$ and $U_j$, for $1\leq i<j\leq \ell$ into $3$-stars such that the resulting 3-star system $(\hat V,\hat {\mathcal B})$ is $k$-chromatic. To do so, we will decompose the edges between $V$ and $U$ into 3-stars in a way such that no 3-subset in $V$ is monochromatic in any putative $(k-1)$-colouring. Three cases arise.

\textbf{Case 1:} $2k-1 \equiv 0$ (mod 3). Then $2k-1=3t$ for some $t \in \mathbb Z^+$. The number of 3-subsets of the set $V$ is $\binom{2k-1}{3}=t(k-1)(2k-3)$. We now fix $\ell=(k-1)(2k-3)$. Partition the set of all 3-subsets of $V$ into $\ell$ sets $\mathbb{T}_1, \ldots, \mathbb{T}_{\ell}$ each consisting of $t$ mutually disjoint 3-subsets. This partition is known to exist by \cite{Baranyai1975}. Let $\mathbb{T}_i=\big\{ \{x_1^i,y_1^i,z_1^i\}, \ldots, \{x_t^i,y_t^i,z_t^i\} \big\}$, $1\leq i \leq \ell$. Decompose the edges between $V$ and $U_i$ into the 3-stars of the set
$\mathcal T_i=\bigcup\limits_{u \in U_i}\big\{ \{u;x_1^i,y_1^i,z_1^i\}, \ldots, \{u;x_t^i,y_t^i,z_t^i\} \big\}$, where $1\leq i \leq \ell$. 

We now begin to decompose the edges between $U_i$ and $U_j$ for $1\leq i<j\leq \ell$ into $3$-stars. If $|C_1|$ is even, then let $r=\frac{|C_1|}{2}$. Observe that $2r=|C_1| \leq |C_2|$ and so for each $p\in\{1,2,\ldots,|C_1|\}$ we decompose the edges between $c_p^1\times \{i\}$ and $U_j$ into the $3$-stars of the set:
$\mathcal S_p^1=\big\{\{c_p^1\times\{i\};c_1^1\times \{j\},c_2^1\times \{j\},c_1^2\times \{j\}\},$ $\{c_p^1\times\{i\};c_3^1\times\{j\},c_4^1\times \{j\},c_2^2\times \{j\}\},$ $\ldots,$ $\{c_p^1\times \{i\};c_{|C_1|-1}^1\times \{j\},c_{|C_1|}^1\times \{j\},c_r^2\times \{j\}\}\big\}$ $\cup$ 
 $\big\{ \{c_p^1\times \{i\};x_t,y_t,z_t\} \mid 1\leq t\leq \frac{n_{k-1}-3r}{3} \big\}$, where the sets $\{x_t,y_t,z_t\}$, $1\leq t\leq \frac{n_{k-1}-3r}{3}$, form a partition of $U_j\setminus \big((C_1\times\{j\}) \cup \{c_1^2\times \{j\},c_2^2\times \{j\},\ldots,c_r^2\times \{j\}\}\big)$.
 If $|C_1|$ is odd, then let $r=\frac{|C_1|+1}{2}$ and decompose the edges between $c_p^1\times \{i\}$ and $U_j$ into the $3$-stars of the set:
$\mathcal S_p^1=\big\{\{c_p^1\times \{i\};c_1^1\times \{j\},c_2^1\times \{j\},c_1^2\times \{j\}\},$ $\{c_p^1\times \{i\};c_3^1\times \{j\},c_4^1\times \{j\},c_2^2\times \{j\}\},$ $\ldots,$ $\{c_p^1\times \{i\};c_{|C_1|-2}^1\times \{j\},c_{|C_1|-1}^1\times \{j\},c_{r-1}^2\times \{j\}\}$, $\{c_p^1\times \{i\};c_{|C_1|}^1\times \{j\},c_r^2\times \{j\},c_{r+1}^2\times \{j\}\}\big\}$ $\cup$  $\big\{ \{c_p^1\times \{i\};x_t,y_t,z_t\} \mid 1\leq t\leq \frac{n_{k-1}-3r}{3} \big\}$,  where the sets $\{x_t,y_t,z_t\}$, $1\leq t\leq \frac{n_{k-1}-3r}{3}$,  form a partition of $U_j\setminus \big((C_1\times\{j\}) \cup \{c_1^2\times \{j\},c_2^2\times \{j\},\ldots,c_{r+1}^2\times \{j\}\}\big)$. The edges between $(C_2\times \{i\}) \cup \cdots \cup (C_{k-2}\times \{i\})$ and $U_j$ are decomposed into sets of $3$-stars $ \mathcal S_1^2,\ldots,\mathcal S_{|C_2|}^2, \ldots, \mathcal S_1^{k-2},\ldots,\mathcal S_{|C_{k-2}|}^{k-2}$ in a similar manner.

Next, we decompose the edges between $C_{k-1}\times \{i\}$ and $U_j$ into $3$-stars. If $|C_{k-1}|=2n$ for some $n \in \mathbb Z^+$, then we let $C_k=\{c_{n+1}^{k-1},\ldots, c_{|C_{k-1}|}^{k-1} \}$. If $2n$ is a multiple of 3, then for each $p\in\{1,2,\ldots,|C_{k-1}|\}$ we decompose the edges between $c_p^{k-1}$ and $U_j$ into the $3$-stars of the set: 
$\mathcal S_p^{k-1}=\big\{ \{c_p^{k-1}\times \{i\};c_1^{k-1}\times \{j\},c_2^{k-1}\times \{j\},c_{n+1}^{k-1}\times \{j\}\}$, $\{c_p^{k-1}\times \{i\};c_{n+2}^{k-1}\times \{j\},c_{n+3}^{k-1}\times \{j\},c_3^{k-1}\times \{j\}\}$, $\ldots,$ $\{c_p^{k-1}\times \{i\};c_{n-2}^{k-1}\times \{j\},c_{n-1}^{k-1}\times \{j\},c_{|C_{k-1}|-2}^{k-1}\times \{j\}\},$ $\{c_p^{k-1}\times \{i\};c_{|C_{k-1}|-1}^{k-1}\times \{j\},c_{|C_{k-1}|}^{k-1}\times \{j\},c_n^{k-1}\times \{j\}\}\big\}$ $\cup$ $\big\{ \{c_p^{k-1}\times \{i\};x_t,y_t,z_t\} \mid 1 \leq t\leq \frac{n_{k-1}-|C_{k-1}|}{3} \big\}$, where  the sets $\{x_t,y_t,z_t\}$, $1 \leq t\leq \frac{n_{k-1}-|C_{k-1}|}{3}$, form a partition of $\bigcup\limits_{s=1}^{k-2}(C_s\times \{j\})$. If $2n$  is congruent to 1 modulo 3, then the edges are decomposed into the 3-stars of the set:
$\mathcal S_p^{k-1}=\big\{ \{c_p^{k-1}\times \{i\};c_1^{k-1}\times \{j\},c_2^{k-1}\times \{j\},c_{n+1}^{k-1}\times \{j\}\}$, $\{c_p^{k-1}\times \{i\};c_{n+2}^{k-1}\times \{j\},c_{n+3}^{k-1}\times \{j\},c_3^{k-1}\times \{j\}\}$, $\ldots,$ $\{c_p^{k-1}\times \{i\};c_{n-2}^{k-1}\times \{j\},c_{n-1}^{k-1}\times \{j\},c_{|C_{k-1}|-3}^{k-1}\times \{j\}\},$ 
$\{c_p^{k-1}\times \{i\};c_{|C_{k-1}|-2}^{k-1}\times \{j\},c_{|C_{k-1}|-1}^{k-1}\times \{j\},c_n^{k-1}\times \{j\}\},$ $\{c_p^{k-1}\times \{i\};c_{|C_{k-1}|}^{k-1}\times \{j\},a,b\}\big\}$ $\cup$ $\big\{ \{c_p^{k-1}\times \{i\};x_t,y_t,z_t\} \mid 1 \leq t\leq \frac{n_{k-1}-(|C_{k-1}|+2)}{3} \big\}$, where $a$ and $b$ are distinct vertices of  $\bigcup\limits_{s=1}^{k-2}(C_s\times \{j\})$ and the sets $\{x_t,y_t,z_t\}$, $1 \leq t\leq \frac{n_{k-1}-(|C_{k-1}|+2)}{3}$, form a partition of $\bigcup\limits_{s=1}^{k-2}(C_s\times \{j\})$ $\setminus \{a,b\}$. The edges are decomposed into 3-stars in a similar manner when $2n$ is congruent to 2 modulo 3.
If $|C_{k-1}|=2n+1$ for some $n \in \mathbb Z^+$, then we let $C_k=\{c_{n+1}^{k-1}\,\ldots, c_{|C_{k-1}|}^{k-1}\}$ and decompose the edges into $3$-stars in a similar manner. 
Now $\mathcal U_j^i=(\bigcup\limits_{a=1}^{|C_1|}\mathcal S_a^1)\cup\cdots\cup(\bigcup\limits_{a=1}^{|C_{k-1}|}\mathcal S_a^{k-1})$ 
forms a 3-star decomposition of the edges between $U_i$ and $U_j$, $1\leq i<j\leq \ell$.

Let $\hat{\mathcal B}=\mathcal B \cup (\bigcup \limits_{i=1}^{\ell} \mathcal A_i)  \cup (\bigcup \limits_{i=1}^{\ell} \mathcal T_i) \cup (\bigcup \limits_{j=2}^{\ell} \mathcal U_j^1) \cup (\bigcup \limits_{j=3}^{\ell} \mathcal U_j^2) \cup \cdots \cup \mathcal U_{\ell}^{\ell -1}$. The $k$ colour classes $\{1,2\} \cup \bigcup\limits_{i=1}^{\ell}(C_1\times\{i\})$, $\{3,4\} \cup \bigcup\limits_{i=1}^{\ell}(C_2\times\{i\})$, $\cdots$, $\{2k-1,2k-2\} \cup \bigcup\limits_{i=1}^{\ell}(C_{k-1}\setminus C_k)\times\{i\})$, $\{2k-1\} \cup \bigcup\limits_{i=1}^{\ell}(C_k\times\{i\})$ exhibit a $k$-colouring of $(\hat V,\hat{\mathcal B})$. Since $(U_i,\mathcal A_i)$, $1\leq i \leq \ell$ is a $(k-1)$-chromatic 3-star system, there are vertices $u_1^i, \ldots, u_{k-1}^i$ in $U_i$ such that $u_s^i\in C_s\times \{i\}$ for each $1\leq s\leq k -1$. Since $|V|=2k-1>2(k-1)$ then when attempting to colour $V$ with $k-1$ colours, some colour must occur at least thrice within $V$. Thus, for some $j\in\{1,\ldots,t\}$ and some $i \in \{1,\ldots,\ell\}$ and for some 
$s\in \{1,\ldots,k-1\}$, the three vertices $x_j^i$, $y_j^i$ and $z_j^i$ are each coloured with the same colour $s$. Then the 3-star $\{u_s^i;x_j^i,y_j^i,z_j^i\}$ would be monochromatic and hence $(\hat V,\hat {\mathcal B})$ cannot be coloured with $k-1$ colours.

Therefore, $(\hat V,\hat{\mathcal B})$ is a $k$-chromatic 3-star system of order $n_k=(2k-1)+n_{k-1}(k-1)(2k-3)$.

\textbf {Case 2:} $2k-1 \equiv 1$ (mod 3). Then $2k-1=3t+1$ for some $t \in \mathbb Z^+$. The number of 3-subsets of the set $V$ is $\binom{2k-1}{3}=\frac{t}{2}(2k-1)(2k-3)$. Note that since $\binom{2k-1}{3}=\frac{t}{2}(2k-1)(2k-3)$ is an integer and $2k-1$ and $2k-3$ are not divisible by $2$, $t$ must be even. Let $t=2t'$.
Partition the set of all 3-subsets of $V$ into $\ell'=(2k-1)(k-2)$ sets of disjoint 3-subsets $\mathbb{T}_1,\ldots,\mathbb{T}_{\ell'}$ of size $t$ and $\ell''=2k-1$ sets of disjoint 3-subsets $\mathbb{T}_{\ell'+1},\ldots,\mathbb{T}_{\ell'+\ell''}$ of size $t'$. Such a partition is known to exist by Theorem~1 in~\cite{Baranyai1975}. We now fix $\ell=\ell'+\ell''$.
Let $\mathbb{T}_i=\big\{ \{x_1^i,y_1^i,z_1^i\}, \ldots, \{x_t^i,y_t^i,z_t^i\} \big\}$, $1\leq i \leq \ell'$ and $\mathbb{T}_{i}=\big\{ \{x_1^{i},y_1^{i},z_1^{i}\}, \ldots, \{x_{t'}^{i},y_{t'}^{i},z_{t'}^{i}\} \big\}$, $\ell'+1\leq i \leq \ell'+\ell''$. 

Now, we decompose the edges between $V$ and $U_i$, $1\leq i\leq \ell'$, into 3-stars.  Let $v$ be the single vertex of the set $V\setminus \bigcup\limits_{j=1}^{t}\{x_j^i,y_j^i,z_j^i\}$. The edges between $v$ and $U_i$, $1\leq i\leq \ell'$, are decomposed into 3-stars in a manner similar to the decomposition of the edges between a vertex $u\in U_i$ and $U_j$, $1\leq i< j\leq \ell$, in Case~1. Let $\mathcal T _i^1$ be the decomposition of the edges between $v$ and $U_i$ into 3-stars, $\mathcal T_i^2=\bigcup\limits_{u\in U_i}\big\{ \{u;x_1^i,y_1^i,z_1^i\}, \ldots, \{u;x_t^i,y_t^i,z_t^i\} \big\}$ and $\mathcal T_i=\mathcal T_i^1 \cup \mathcal T_i^2$, $1\leq i \leq \ell'$. Then $\mathcal T_i$ is a decomposition of the edges between $V$ and $U_i$, $1\leq i\leq \ell'$, into 3-stars. The edges between $V$ and $U_{i}$ are decomposed into 3-stars $\mathcal T_{i}$ in a similar manner for $\ell'+1\leq i \leq \ell'+\ell''$. 

The edges between $U_i$ and $U_j$, $1\leq i< j\leq \ell$, are decomposed into sets of 3-stars $\mathcal U_j^i$ in a manner similar to Case~1.

Let $\hat{\mathcal B}=\mathcal B \cup (\bigcup \limits_{i=1}^{\ell} \mathcal A_i)  \cup (\bigcup \limits_{i=1}^{\ell} \mathcal T_i) \cup (\bigcup \limits_{j=2}^{\ell} \mathcal U_j^1) \cup (\bigcup \limits_{j=3}^{\ell} \mathcal U_j^2) \cup \cdots \cup \mathcal U_{\ell}^{\ell -1}$. The $k$ colour classes $\{1,2\} \cup \bigcup\limits_{i=1}^{\ell}(C_1\times\{i\})$, $\{3,4\} \cup \bigcup\limits_{i=1}^{\ell}(C_2\times\{i\})$, $\cdots$, $\{2k-1,2k-2\} \cup \bigcup\limits_{i=1}^{\ell}(C_{k-1}\setminus C_k)\times\{i\})$, $\{2k-1\} \cup \bigcup\limits_{i=1}^{\ell}(C_k\times\{i\})$ exhibit a $k$-colouring of $(\hat V,\hat{\mathcal B})$. It is seen that $(\hat V,\hat{\mathcal B})$ is a $k$-chromatic 3-star system of order $n_k=(2k-1)+n_{k-1}(2k-1)(k-1)$ in a manner similar to Case~1.

\textbf {Case 3:} $2k-1 \equiv 2$ (mod 3). Then $2k=3t$ for some $t \in \mathbb Z^+$. The number of 3-subsets of the set $V$ is $\binom{2k}{3}=t(2k-1)(k-1)$. We now fix $\ell=(2k-1)(k-1)$. Partition the set of all 3-subsets of $V$ into $\ell$ sets of disjoint 3-subsets $\mathbb{T}_1, \ldots, \mathbb{T}_{\ell}$ of size $t$. Let $\mathbb{T}_i=\big\{ \{x_1^i,y_1^i,z_1^i\}, \ldots, \{x_t^i,y_t^i,z_t^i\} \big\}$, $1\leq i \leq \ell$. We decompose the edges between $V$ and $U_i$ into the 3-stars of the set
$\mathcal T_i=\bigcup\limits_{u\in U_i}$$\big\{ \{u;x_1^i,y_1^i,z_1^i\},$ $\ldots,$ $\{u;x_t^i,y_t^i,z_t^i\} \big\}$, where $1\leq i \leq \ell$. The edges between $U_i$ and $U_j$, $1\leq i< j\leq \ell$, are decomposed into 3-stars $\mathcal U_j^i$ in a manner similar to Case~1. Let $\hat{\mathcal B}=\mathcal B \cup (\bigcup \limits_{i=1}^{\ell} \mathcal A_i)  \cup (\bigcup \limits_{i=1}^{\ell} \mathcal T_i) \cup (\bigcup \limits_{j=2}^{\ell} \mathcal U_j^1) \cup (\bigcup \limits_{j=3}^{\ell} \mathcal U_j^2) \cup \cdots \cup \mathcal U_{\ell}^{\ell -1}$. The $k$ colour classes $\{1,2\} \cup \bigcup\limits_{i=1}^{\ell}(C_1\times\{i\})$, $\{3,4\} \cup \bigcup\limits_{i=1}^{\ell}(C_2\times\{i\})$, $\cdots$, $\{2k-1,2k-2\} \cup \bigcup\limits_{i=1}^{\ell}(C_{k-1}\setminus C_k)\times\{i\})$, $\{2k-1,2k\} \cup \bigcup\limits_{i=1}^{\ell}(C_k\times\{i\})$ exhibit a $k$-colouring of $(\hat V,\hat{\mathcal B})$. It is seen that $(\hat V,\hat{\mathcal B})$ is a $k$-chromatic 3-star system of order $n_k=2k+n_{k-1}(2k-1)(k-1)$ in a manner similar to Case~1.
\end{Proof}

We finish this section with the following corollary.
\begin{Corollary}\label{cor2.1}
For any integer $k \geq 2$, there exists some integer $n_k$ such that for all admissible $n \geq n_k$, there exists a $k$-chromatic $3$-star system of order $n$.
\end{Corollary}

\begin{Proof}
For $k=2$, apply Theorem~\ref{thm2}. For $k\geq 3$, apply Theorem~\ref{thm1} and Theorem~\ref{thm3}.
\end{Proof}

Observe that the cardinalities of the colour classes of the $k$-chromatic 3-star systems from Theorem~\ref{thm1} and Theorem~\ref{thm3} are not excessively imbalanced. The cardinalities of the colour classes of the $k$-chromatic 3-star systems either differ by at most a small number or the cardinality of the largest colour class of the $k$-chromatic 3-star systems is almost two times that of the cardinality of the smallest colour class. 

\section{$\mathbf k$-colourings of $\mathbf e$-star systems}
In this section, we generalize the results of the last section for $e$-star systems for all $e\geq 3$. For any arbitrary integers $e\geq 3$ and $k \geq 2$, we show that there exists some integer $n_k$ where $n_k \equiv 0$ (mod $2e$) such that for all $n \geq n_k$ where $n\equiv 0,1$ (mod $2e$), there exists a $k$-chromatic $e$-star system of order $n$. We first construct a strongly equitable 2-chromatic $e$-star system of order $2e$ for all $e\geq 3$. 

\begin{Theorem}\label{thm4}
There exists a strongly equitable 2-chromatic $e$-star system of order $2e$ for all $e\geq 3$.
\end{Theorem}

\begin{Proof}
By Theorem~\ref{thm2}, there exists an equitably 2-chromatic 3-star system of order six, $(V,\mathcal B)$, where $V=\{1,2,3,4,5,6\}$ and $\mathcal B=\big \{ \{1;3,5,6\}$, $\{2;1,3,6\}$, $\{4;1,2,3\}$, $\{5;2,3,4\}$, $\{6;3,4,5\}\big\}$ with colour classes $R=\{1,3,5\}$ and $Y=\{2,4,6\}$. We construct a 2-chromatic 4-star system of order eight, $(\hat{V},\hat{\mathcal B})$, from $(V,\mathcal B)$. Let $\hat{V}=V\cup \{7,8\}$ and $\hat{\mathcal B}=$$\big\{ \{1;3,5,6,8\}$, $\{2;1,3,6,8\}$, $\{4;1,2,3,8\}$, $\{5;2,3,4,7\}$, $\{6;3,4,5,7\}$, $\{7;1,2,3,4\}$,  $\{8;3,5,6,7\}\big\}$, $\hat{R}=\{1,3,5,7\}$ and $\hat{Y}=\{2,4,6,8\}$. Then $(\hat{V},\hat{\mathcal B})$ is an equitably 2-chromatic 4-star system of order 8 with colour classes $\hat{R}$ and $\hat{Y}$.

We now generalize this construction, which can be used in an iterative manner. Suppose that there exists an equitably 2-chromatic $e$-star system of order $2e$, $(V,\mathcal B)$, where $V=\{1,\ldots,2e\}$ is the set of points which is partitioned into two subsets $R$ and $Y$ and $\mathcal B=\big\{ \{1;a_1^1,\ldots,a_e^1\},$ $\{2;a_1^2,\ldots,a_e^2\},$ $\{4;a_1^4,\ldots,a_e^4\},$ $\{5;a_1^5,\ldots,a_e^5\},$ $\ldots, \{2e;a_1^{2e},\ldots,a_e^{2e}\} \big\}$ is the set of blocks. We construct an equitably 2-chromatic $(e+1)$-star system of order $2e+2$, $(\hat{V},\hat{\mathcal B})$, from $(V,\mathcal B)$. Let $\hat{V}=V\cup \{2e+1,2e+2\}$ and $\hat{\mathcal B}=\big\{ \{1;a_1^1,\ldots,a_e^1,2e+2\},$ $\{2;a_1^2,\ldots,a_e^2,2e+2\},$ $\{4;a_1^4,\ldots,a_e^4,2e+2\},$ $\{5;a_1^5,\ldots,a_e^5,2e+2\},$ $\ldots, \{e+1;a_1^{e+1},\ldots,a_e^{e+1},2e+2\}$, $\{e+2;a_1^{e+2},\ldots,a_e^{e+2},2e+1\}$, $\{e+3;a_1^{e+3},\ldots,a_e^{e+3},2e+1\},$ $\ldots, \{2e;a_1^{2e},\ldots,a_e^{2e},2e+1\}$, $\{2e+1;1,2,\ldots,e+1\}$, $\{2e+2;3,e+2,e+3,\ldots,2e,2e+1\} \big\}$. Then $(\hat{V},\hat{\mathcal B})$ is a strongly equitable 2-chromatic $(e+1)$-star system of order $2e+2$ with colour classes $\hat{R}=\{1,3,5,\ldots,2e+1\}$ and $\hat{Y}=\{2,4,6,\ldots,2e+2\}$.
\end{Proof}
We now show how to construct a $k$-chromatic $e$-star system from a smaller $k$-chromatic $e$-star system.
\begin{Theorem}\label{thm5}
Let $k\geq 2$ and $e\geq 3$. If there exists a $k$-chromatic $e$-star system of order $n_0$ such that $n_0\equiv 0,1$ (mod $2e$), then there exists a $k$-chromatic $e$-star system for all $n >n_0$ such that $n\equiv 0,1$ (mod $2e$).
\end{Theorem}
\begin{Proof}
Suppose that there exists a $k$-chromatic $e$-star system, $(V,\mathcal B)$, where $V=\{1,\ldots, n_0\}$ is the set of points and $\mathcal B$ is the set of blocks and $n_0\equiv 0$ or 1 (mod $2e$). 
 Given a $k$-colouring of $(V,\mathcal B)$, with colours $1,2,\ldots,k$ and colour classes $C_1,C_2,\ldots,C_k$, let $R=\bigcup\limits_{i=1}^{\ell}C_i$ and $Y=\bigcup\limits_{i=\ell+1}^{k}C_i$ for some integer $\ell$ such that $1\leq \ell < k$. 
Let $r_1, r_2,\ldots,r_{|R|}$ be the elements of $R$ and $y_1, y_2, \ldots, y_{|Y|}$ be the elements of $Y$. Observe that $R$ and $Y$ form a partition of $V$. Without loss of generality assume that $|R|\geq |Y|$. 

\textbf {Case 1.} Suppose that $n_0\equiv 0$ (mod $2e$). Then $n_0=2et$, $t\geq 1$.

 First, we construct a $k$-chromatic $S_e(2et+1)$, $(\hat{V}, \hat{\mathcal B})$, from $(V,\mathcal B)$. Let $\hat{V}=V\cup \{2et+1\}$.

If $|R|=|Y|$, let
$\hat{\mathcal B}=\mathcal B \cup \mathcal T$ where $\mathcal T$ is the set
$\big \{\{2et+1;r_1,\ldots,r_{e-1},y_1\}$, $\{2et+1;y_2,\ldots,y_e,r_e\}$, $\ldots,$ $\{2et+1;r_{|R|-(e-1)},\ldots,r_{|R|-1},y_{|Y|-(e-1)}\},$ $\{2et+1;y_{|Y|-e+2},\ldots,y_{|Y|},r_{|R|}\}\big \}.$

Otherwise, $|R|>|Y|$. Then let
$\hat{\mathcal B}=\mathcal B \cup \mathcal T$ where $\mathcal T$ is the set
$\big \{\{2et+1;y_1,\ldots,y_{e-1},r_1\},$ $\ldots,$ $\{2et+1;y_{t'e-(t'-1)},\ldots,y_{(t'+1)e-t'-1},r_{t'+1}\},$ $\{2et+1;y_{(t'+1)e-t'},\ldots,y_{|Y|},r_{t'+2},\ldots,r_{t'+i-1}\}$, $\{2et+1;r_{t'+i},\ldots,r_{t'+i+e-1}\}$, $\ldots,$ $\{2et+1;r_{|R|-(e-1)},\ldots,r_{|R|}\}\big\},$ where $t'=\lfloor\frac{|Y|}{e-1}\rfloor-1$ and $i=(t'+2)e-|Y|-t'+1$. 

 Note that $(\hat{V}, \hat{\mathcal B})$ is not $(k-1)$-colourable because it contains a copy of $(V,\mathcal B)$. Observe that $C_1,\ldots,C_{k-1},C_k\cup \{2et+1\}$ are the colour classes of a $k$-colouring of $(\hat{V},\hat {\mathcal B})$. Therefore, $(\hat{V},\hat {\mathcal B})$ is a $k$-chromatic $e$-star system of order $2et+1$.

Next, we construct a $k$-chromatic $S_e(2et+2e)$, $(\tilde{V}, \tilde{\mathcal B})$, from $(V,\mathcal B)$. Let $\tilde{V}=V\cup V'$ where $V'=\{2et+1,\ldots,2et+2e\}$ and let $(V',\mathcal B')$ be a 2-chromatic $e$-star system of order $2e$ constructed from Theorem~\ref{thm4} with colour classes $R'=\{2et+1,2et+3,\ldots,2et+2e-1\}$ and $Y'=\{2et+2,2et+4,\ldots,2et+2e\}$. Let $\tilde{\mathcal B}=\mathcal B \cup \mathcal B' \cup \mathcal T$ where $\mathcal T=\bigcup\limits_{v\in V}\big\{\{v;2et+1,\ldots,2et+e\},\{v;2et+e+1,\ldots,2et+2e\}\big\}$. Note that $(\tilde{V}, \tilde{\mathcal B})$ is not $(k-1)$-colourable because it contains a copy of $(V,\mathcal B)$. Observe that $C_1\cup R',C_2\cup Y',C_3,\ldots,C_k$ are the colour classes of a $k$-colouring of $(\tilde{V},\tilde {\mathcal B})$. Therefore, $(\tilde{V},\tilde {\mathcal B})$ is a $k$-chromatic $e$-star system of order $2et+2e$.

\textbf {Case 2.} Suppose that $n_0\equiv 1$ (mod $2e$). Then $n_0=2et+1$, $t\geq 1$.
We construct a $k$-chromatic $S_e(2et+2e)$, $(\hat{V}, \hat{\mathcal B})$, from $(V,\mathcal B)$. Let $\hat{V}=V\cup V_1\cup V_2$ where $V_1=\{2et+2,\ldots,2et+e+1\}$ and $V_2=\{2et+e+2,\ldots,2et+2e\}$. Let $v_0=r_{|R|}$ and $R'=R\setminus \{v_0\}$.

If $|R'|=|Y|$, let $\mathcal T_1=\bigcup\limits_{v\in V_2}\big\{\{v;r_1,\ldots,r_{e-1},y_1\},$ $\{v;y_2,\ldots,y_e,r_e\},\ldots,\{v;r_{|R|-e},$ $\ldots,$ $r_{|R|-2},$ $y_{|Y|-e+1}\},$ $\{v;y_{|Y|-e+2},\ldots,y_{|Y|},r_{|R|-1}\}\big\}$.

If $|R'|=|Y|-1$ or $|R'|>|Y|$, decompose the edges between $V_2$ and $V\setminus \{v_0\}$ into a set of $e$-stars $\mathcal T_1$ in a manner similar to the case $|R'|=|Y|$.

Let $V'=V_1\cup V_2 \cup \{v_0\}$ and $(V',\mathcal B')$ be a 2-chromatic $e$-star system of order $2e$ constructed from Theorem~\ref{thm4} with colour classes $R'=\{2et+3,2et+4,\ldots,2et+e+1\}$ and $Y'=\{2et+2,2et+e+2,2et+e+3,\ldots,2et+2e\}$. Let $\mathcal T_2=\bigcup\limits_{v\in V\setminus \{v_0\}}\big\{\{v;2et+2,\ldots,2et+e+1\}\big\}$. Let $\hat{\mathcal B}=\mathcal B\cup \mathcal B'\cup \mathcal T_1\cup \mathcal T_2$. Note that $(\hat{V}, \hat{\mathcal B})$ is not $(k-1)$-colourable because it contains a copy of $(V,\mathcal B)$. Observe that $C_1\cup R',C_2,\ldots, C_{k-1},C_k\cup Y'$ are the colour classes of a $k$-colouring of $(\hat{V},\hat {\mathcal B})$. Therefore, $(\hat{V},\hat {\mathcal B})$ is a $k$-chromatic $e$-star system of order $2et+2e$.
\end{Proof}
We now show how to iteratively construct a $k$-chromatic $e$-star system from a $(k-1)$-chromatic $e$-star system.
\begin{Theorem}\label{thm6}
	Let $k\geq 3$ and $e\geq 3$. If there exists a $(k-1)$-chromatic $e$-star system of order $n_{k-1}\equiv 0$ (mod $2e$) then there exists a $k$-chromatic $e$-star system of order $n_k \equiv 0$ (mod $2e$).
\end{Theorem}
\begin{Proof}
	Let $(U_0,\mathcal A_0)$ be a $(k-1)$-chromatic $e$-star system of order $n_{k-1}$ with colour classes $C_1,\ldots,C_{k-1}$ such that $n_{k-1}\equiv 0$ (mod $2e$), $|C_1| \leq \cdots \leq |C_{k-1}|$ and each vertex of $C_s$ has colour $s$, for $1\leq s\leq k-1$. For each $s\in \{1,\ldots,k-1\}$, let $C_s=\{c_1^s,\ldots,c_{|C_s|}^s\}$. For a positive integer $\ell$ that will be fixed later and for each $i\in \{1,\ldots,\ell\}$, let $U_i=U_0 \times \{i\}$ and $\mathcal A_i = \mathcal A_0\times \{i\}$ where $\mathcal A_0\times \{i\}$ denotes the set $\{S\times \{i\} \mid S\in \mathcal A_0\}$. So $(U_i,\mathcal A_i)$ has colour classes $C_1\times\{i\},\ldots,C_{k-1}\times\{i\}$. Let $W=\{w_1,w_2,\ldots,w_{(e-1)(k-1)+1}\}$, $B=\{b_1,b_2,\ldots,b_{k-2}\}$ and $D=\{d_1,d_2,\ldots,d_e\}$ be pairwise disjoint sets. If $k-1$ is even then let $V=W\cup B$; otherwise let $V=W\cup B\cup D$. Also let $U=\bigcup\limits_{i=1}^{\ell} U_i$ such that $V\cap U=\emptyset$. We will embed $(U_1,\mathcal A_1),\ldots,(U_{\ell},\mathcal A_{\ell})$ into a $k$-chromatic $e$-star system $(\hat V,\hat {\mathcal B})$ where $\hat V=V \cup U$. Let $n_k=|\hat{V}|$ and observe that $n_k\equiv 0$ (mod $2e$). By Theorem~\ref{thm4} and Theorem~\ref{thm5}, the edges of the complete graph on the set $V$ admit a decomposition into $e$-stars which is 2-chromatic. Let $(V,\mathcal B)$ be a 2-chromatic $e$-star system on the set $V$ (so ($V,\mathcal B)$ will not obstruct the $k$-colouring that we will construct). We now need to decompose the edges between $V$ and $U$ and the edges between $U_i$ and $U_j$, for $1\leq i<j\leq \ell$ into $e$-stars such that the resulting $e$-star system $(\hat V,\hat {\mathcal B})$ is $k$-chromatic. To do so, we will decompose the edges between $V$ and $U$ into $e$-stars in a way such that no $e$-subset in $V$ is monochromatic in any putative $(k-1)$-colouring. 
	
	The number of $e$-subsets of the set $W$ is $N=\binom{(e-1)(k-1)+1}{e}$. Let $a=(k-2)+ \lfloor \frac{-k+2}{e}\rfloor$ and $N=aq+r$ where $q$ is a positive integer and $0\leq r<a$ is an integer. Partition the set of all $e$-subsets of $W$ into $q$ sets of disjoint $e$-subsets $\mathbb{T}_1,\ldots,\mathbb{T}_q$ of size $a$ and one set of disjoint $e$-subsets $\mathbb{T}_{q+1}$ of size $r$. Such a partition is known to exist by Theorem~1 in~\cite{Baranyai1975}. We now fix $\ell=q+1$ if $r>0$ and $\ell=q$ if $r=0$. Let $\mathbb{T}_i=\big\{\{(x_{1})_1^i,\ldots,(x_{e})_1^{i}\},\ldots,\{(x_{1})_a^i,\ldots,(x_{e})_a^{i}\}\big\}$, $1\leq i\leq q$ and $\mathbb{T}_{q+1}=\big\{\{(x_{1})_1^{q+1},\ldots,(x_{e})_1^{q+1}\},\ldots,\{(x_{1})_r^{q+1},\ldots,(x_{e})_r^{q+1}\}\big\}$ if $r>0$ and $\mathbb{T}_{q+1}=\emptyset$ if $r=0$. Let $W_i=W\setminus \bigcup\limits_{j=1}^{a}\big(\{(x_{1})_j^{i},\ldots,(x_{e})_j^{i}\}\big)$ for $1\leq i\leq q$ and $W_{q+1}=W\setminus \bigcup\limits_{j=1}^{r}\big(\{(x_{1})_j^{q+1},\ldots,(x_{e})_j^{q+1}\}\big)$ if $r>0$ and $W_{q+1}=\emptyset$ if $r=0$.

	For each $i\in \{1,2,\ldots,\ell\}$ and for each $w\in W_i$, we decompose the edges between $w$ and $U_i$. First, we decompose all the edges between $w$ and $C_1\times\{i\}$ into $e$-stars. Let $\mathcal R_{w}^{1,2}=\big\{\{w;c_1^1\times \{i\},\ldots,c_{e-1}^{1}\times \{i\},c_1^2\times \{i\}\},\ldots,
	\{w;c_{(t-1)(e-1)+1}^1\times \{i\},\ldots,c_{t(e-1)}^{1}\times \{i\},c_{t}^2\times \{i\}\},\{w;c_{t(e-1)+1}^1\times\{i\},\ldots,c_{|C_1|}^1\times \{i\},c_{t+1}^2\times \{i\},\ldots,c_{t+t'}^2\times \{i\}\}\big\}$ where $t$ and $t'$ are integers such that  $t=\lfloor \frac{|C_1|}{e-1}\rfloor$ and $t'=(t+1)e-|C_1|$. Then $\mathcal R_{w}^{1,2}$ is a decomposition of all the edges between $w$ and $C_1\times \{i\}$ along with the edges between $w$ and $C'_2\times \{i\}$ into $e$-stars where $C'_2=\big\{c_{1}^2,\ldots,c_{t+t'}^2\big\}\subseteq C_2$. 
	For each $j=2,3,\ldots,k-2$, iteratively proceed in a similar manner to decompose all the edges between $w$ and $(C_j\setminus C'_j)\times \{i\}$ into a set of $e$-stars $\mathcal R_w^{j,j+1}$. So $\bigcup\limits_{j=1}^{k-2} \mathcal R_{w}^{j,j+1}$ gives a decomposition of the edges between $w$ and $(C_1\cup \cdots \cup C_{k-2}\cup C'_{k-1})\times \{i\}$ where $C'_{k-1}\subseteq C_{k-1}$. Let $m=\frac{|C_{k-1}\setminus C'_{k-1}|}{e}$ and let $C''_{k-1}=C_{k-1}\setminus C'_{k-1}$ and without loss of generality we can let $C''_{k-1}=\{c_1^{k-1},\ldots,c_{me}^{k-1}\}$. Without loss of generality let $C_k=\big\{c_1^{k-1},c_{e+1}^{k-1},c_{2e+1}^{k-1},\ldots,c_{(m-1)e+1}^{k-1}\big\}$ and decompose the edges between $w$ and $C''_{k-1}\times \{i\}$ into $e$-stars of the set: $\mathcal R_w^{k-1}=\big\{\{w;c_1^{k-1}\times \{i\},c_2^{k-1}\times \{i\},\ldots,c_e^{k-1}\times \{i\}\},$ $\{w;c_{e+1}^{k-1}\times \{i\},c_{e+2}^{k-1}\times \{i\},\ldots,c_{2e}^{k-1}\times \{i\}\},$ $\ldots,$ $\{w;c_{(m-1)e+1}^{k-1}\times \{i\},c_{(m-1)e+2}^{k-1}\times \{i\},\ldots,c_{me}^{k-1}\times \{i\}\}\big\}$. 
	
	For each $i\in \{1,2,\ldots,q\}$, let $\mathcal T_i=\bigcup\limits_{u\in U_i}\big\{\{u;(x_{1})_1^i,\ldots,(x_{e})_1^{i}\},\ldots,\{u;(x_{1})_a^i,\ldots,(x_{e})_a^{i}\}\big\}$ and $\mathcal T_{q+1}=\big\{\{u;(x_{1})_1^{q+1},\ldots,(x_{e})_1^{q+1}\},$ $\ldots,$ $\{u;(x_{1})_r^{q+1},\ldots,(x_{e})_r^{q+1}\}\big\}$ if $r>0$ and $\mathcal T_{q+1}=\emptyset$ if $r=0$. Then $\mathcal P_i= \mathcal T_i \cup \bigcup\limits_{w\in W_i}\big(\mathcal R_w^{k-1}\cup (\bigcup\limits_{j=1}^{k-2}\mathcal R_w^{j,j+1})\big)$ is a decomposition of the edges between $W$ and $U_i$ into $e$-stars. The edges between $U_i$ and $U_j$, $1\leq i < j \leq \ell$ are decomposed into sets of $e$-stars $\mathcal U_j^i$ in a manner similar to the decomposition of the edges between $W_i$ and $U_i$, $1\leq i\leq \ell$. 
	
	We have two cases.
	
	\textbf {Case 1.} $k-1$ is even. Then $V=W\cup B$. The edges between $B$ and $U_i$, $1\leq i\leq \ell$ are decomposed into sets of $e$-stars $\mathcal F_i$ in a manner similar to the decomposition of the edges between $W_i$ and $U_i$, $1\leq i\leq \ell$. Therefore, $(\hat V,\hat{\mathcal B})$ where $\hat{\mathcal B}=\mathcal B \cup (\bigcup \limits_{i=1}^{\ell} \mathcal A_i)  \cup (\bigcup\limits_{i=1}^{\ell}\mathcal P_i) \cup (\bigcup\limits_{i=1}^{\ell}\mathcal F_i) \cup (\bigcup \limits_{j=2}^{\ell} \mathcal U_j^1) \cup (\bigcup \limits_{j=3}^{\ell} \mathcal U_j^2) \cup \cdots \cup \mathcal U_{\ell}^{\ell -1}$ is an $e$-star system of order $n_k=(k-1)e+n_{k-1}(q+1)$. A $k$-colouring is given by
	$(\bigcup\limits_{i=1}^{\ell} C_1\times \{i\}) \cup \{w_1,\ldots,w_{e-1},b_1\}$, $(\bigcup\limits_{i=1}^{\ell}C_2\times \{i\}) \cup \{w_e,\ldots,w_{2e-2},b_2\}$, $\ldots, (\bigcup\limits_{i=1}^{\ell}C_{k-2}\times \{i\})\cup \{w_{(e-1)(k-1)-2e-1},\ldots,w_{(e-1)(k-1)-e-1},b_{k-2}\}$, $\big (\bigcup\limits_{i=1}^{\ell}(C_{k-1}\setminus C_k)\times \{i\}\big ) \cup \{w_{(e-1)(k-1)-e},\ldots,w_{(e-1)(k-1)}\}$, $(\bigcup\limits_{i=1}^{\ell}C_k \times \{i\}) \cup \{w_{(e-1)(k-1)+1}\}$.

	\textbf {Case 2.} $k-1$ is odd. then $V=W\cup B\cup D$. The edges between $D$ and $U_i$, $1\leq i\leq \ell$ are decomposed into sets of $e$-stars $\mathcal H_i$ in a manner similar to the decomposition of the edges between $W_i$ and $U_i$, $1\leq i\leq \ell$. Therefore, $(\hat V,\hat{\mathcal B})$ where $\hat{\mathcal B}=\mathcal B \cup (\bigcup \limits_{i=1}^{\ell} \mathcal A_i)  \cup (\bigcup\limits_{i=1}^{\ell}\mathcal P_i) \cup (\bigcup\limits_{i=1}^{\ell}\mathcal F_i) \cup (\bigcup\limits_{i=1}^{\ell}\mathcal H_i) \cup (\bigcup \limits_{j=2}^{\ell} \mathcal U_j^1) \cup (\bigcup \limits_{j=3}^{\ell} \mathcal U_j^2) \cup \cdots \cup \mathcal U_{\ell}^{\ell -1}$ is an $e$-star system of order $n_k=ek+n_{k-1}(q+1)$. A $k$-colouring is given by 
	$(\bigcup\limits_{i=1}^{\ell} C_1\times \{i\})\cup \{w_1,\ldots,w_{e-1},b_1\}$, $(\bigcup\limits_{i=1}^{\ell} C_2\times \{i\}) \cup \{w_e,\ldots,w_{2e-2},b_2\}$ $\ldots, (\bigcup\limits_{i=1}^{\ell} C_{k-2}\times \{i\})\cup \{w_{(e-1)(k-1)-2e-1},\ldots,w_{(e-1)(k-1)-e-1},b_{k-2}\}$, $\big (\bigcup\limits_{i=1}^{\ell}(C_{k-1}\setminus C_k)\times \{i\}\big )\cup \{w_{(e-1)(k-1)-e},\ldots,w_{(e-1)(k-1)},$ $d_1\}$, $(\bigcup\limits_{i=1}^{\ell}C_k \times \{i\})$ $\cup \{w_{(e-1)(k-1)+1},d_2,\ldots,d_e\}$.

Since $(U_i,\mathcal A_i)$, $1\leq i \leq \ell$ is a $(k-1)$-chromatic $e$-star system, there are vertices $u_1^i, \ldots, u_{k-1}^i$ in $U_i$ such that $u_s^i\in C_s\times \{i\}$ for each $1\leq s\leq k -1$. Since $|W|=(e-1)(k-1)+1>(e-1)(k-1)$ then when attempting to colour $W$ with $k-1$ colours, some colour must occur at least $e$ times within $W$. Thus, for some $i \in \{1,\ldots,q\}$ and some $j\in\{1,\ldots,a\}$ and for some $s\in \{1,\ldots,k-1\}$, the vertices $(x_1)^i_j, (x_2)^i_j, \ldots, (x_e)^i_j$ are each coloured with the same colour $s$. Also, for $i=q+1$ and some $j\in\{1,\ldots,r\}$ and for some $s\in \{1,\ldots,k-1\}$, the vertices $(x_1)^i_j, (x_2)^i_j, \ldots, (x_e)^i_j$ are each coloured with the same colour $s$. Then the $e$-star $\{u_s^i;(x_1)^i_j,(x_2)^i_j,\ldots,(x_e)^i_j\}$ would be monochromatic and hence in either case $(\hat V,\hat {\mathcal B})$ cannot be coloured with $k-1$ colours.
\end{Proof}

We summarize this section with the following corollary.
\begin{Corollary}\label{cor3.1}
For any integers $k \geq 2$ and $e\geq 3$, there exists some integer $n_k$ where $n_k \equiv 0$ (mod $2e$) such that for all admissible $n \geq n_k$ where $n\equiv 0,1$ (mod $2e$), there exists a $k$-chromatic $e$-star system of order $n$.
\end{Corollary}

\begin{Proof}
Apply Theorem~\ref{thm4}, Theorem~\ref{thm5} and Theorem~\ref{thm6}.
\end{Proof}

\section{Unique colourings of $\mathbf e$-star systems}
We now investigate uniquely $k$-chromatic $e$-star systems for any $e\geq 3$. We commence by showing that such designs do indeed exist for $k=2$.

\begin{Theorem}\label{thm4.1}
For any $e\geq 3$, there exists a strongly equitable uniquely $2$-chromatic $e$-star system of order $n$ for some $n \equiv 0$ (mod $2e$).
\end{Theorem}
\begin{Proof}
let $(U_0,\mathcal A_0)$ be a strongly equitable 2-chromatic $e$-star system of order $2e$ constructed from Theorem~\ref{thm4} with colour classes $C_1$ and $C_2$ such that each vertex of $C_s$ has colour $s$, for $s=1,2$. For a positive integer $\ell$ that will be fixed later and for each $i\in \{1,\ldots,\ell\}$, let $U_i=U_0\times \{i\}$ and $\mathcal A_i=\mathcal A_0\times \{i\}$ where $\mathcal A_0 \times \{i\}$ denotes the set $\big\{ S\times \{i\} \, | \, S\in \mathcal A_0\big\}$. So $(U_i,\mathcal A_i)$ has colour classes $C_1\times \{i\}$ and $C_2\times \{i\}$. Let $U=\bigcup\limits_{i=1}^{\ell} U_i$.
  Let $A^1=\{a_1^1,\ldots,a_{e}^1\}$, $A^2=\{a_1^2,\ldots,a_{e}^2\}$ and $A=A^1\cup A^2$. Similarly let $F^1=\{f_1^1,\ldots,f_{e}^1\}$, $F^2=\{f_1^2,\ldots,f_{e}^2\}$, $F=F^1\cup F^2$, $G^1=\{g_1^1,\ldots,g_{e}^1\}$, $G^2=\{g_1^2,\ldots,g_{e}^2\}$, $G=G^1\cup G^2$,  $H^1=\{h_1^1,\ldots,h_{e}^1\}$, $H^2=\{h_1^2,\ldots,h_{e}^2\}$, $H=H^1\cup H^2$, $K^1=\{k_1^1,\ldots,k_e^1\}$, $K^2=\{k_1^2,\ldots,k_e^2\}$ and $K=K^1 \cup K^2$. We construct a uniquely 2-chromatic $e$-star system $(\hat{V},\hat{\mathcal B})$ where $\hat{V}=A\cup F\cup G\cup H\cup U\cup K$. 
  
  We need to decompose the edges between the subsets of $\hat{V}$ into $e$-stars such that the resulting $e$-star system $(\hat{V},\hat{\mathcal B})$ is uniquely $2$-chromatic. 
  
  First, we decompose the edges between $A$ and $U$ into $e$-stars in a way such that no $e$-subset in $A$ other than the $e$-subsets $A^1$ and $A^2$ are monochromatic in any putative $2$-colouring of $(\hat{V},\hat{\mathcal B})$. Let $\mathcal D$ be the set of all $e$-subsets of $A$ except 
  $A^1$ and $A^2$. The number of $e$-subsets of the set $A$ is $N=\binom{2e}{e}$ and $|\mathcal D|=N-2$. Let $T_i$, $1\leq i\leq |\mathcal D|$, be the elements of $\mathcal D$. Let $\tilde{A_i}=A\setminus T_i$, $1\leq i\leq |\mathcal D|$. We now fix $\ell=|\mathcal D|$.
  
  Partition $U_0$ into $e$-subsets $E_1$ and $E_2$  such that $E_j \cap C_1 \neq \emptyset$ and $E_j \cap C_2 \neq \emptyset$, $j=1,2$. For each $i\in \{1,2,\ldots,\ell\}$ and for each $w\in \tilde {A_i}$, we decompose the edges between $w$ and $U_i$ into the set of $e$-stars $\mathcal R_w=\bigcup\limits_{j=1}^{2}\big\{ \{w;E_j\times \{i\} \}\big\}$. For each $i\in \{1,\ldots,\ell\}$, let $\mathcal T_i=\bigcup \limits_{u\in U_i} \big\{ \{u;T_i\}\big\}$. Then $\mathcal P_i=\mathcal T_i\cup \big (\bigcup \limits_{w\in \tilde {A_i}} \mathcal R_w \big)$ is a decomposition of the edges between $A$ and $U_i$ and $\bigcup\limits_{i=1}^{\ell} \mathcal P_i$ is a decomposition of the edges between $A$ and $U$ into $e$-stars.
  
  Recall that $(U_i,\mathcal A_i)$ is 2-chromatic. So for any 2-colouring of $(U_i,\mathcal A_i)$ there must be a vertex $u_s^i$ of each colour $s$ in $U_i$, $1\leq s\leq 2$. Note that for each $i\in \{1,\ldots,|\mathcal D|\}$, $\{u_s^i;T_i\}\in \mathcal T_i$ and so to ultimately obtain a valid 2-colouring of $(\hat{V},\hat{\mathcal B})$, the $e$ vertices of $T_i$ must not all have colour $s$, $1\leq s\leq 2$. If a colour $s$ occurs on more than $e$ vertices of $A$ then there exists a $T_i$ with only colour $s$ and then $\{u_s^i;T_i\}$ would be monochromatic. So in any valid 2-colouring of $(\hat{V},\hat{\mathcal B})$, each colour must occur on exactly $e$ points of $A$. If the points of colour $s$ are not all in $A^1$ or $A^2$ then they are the leaves of a monochromatic $e$-star with centre $u_s^i$. Without loss of generality we therefore assign colour $s$ to all $e$ points of $A^s$, $1\leq s\leq 2$.

We now begin to decompose the edges between $A$ and $F$ into $e$-stars. For each $i\in \{1,\ldots,e\}$, let $\mathcal K_i^1=\big\{ \{f_i^1;A^2\}\big\}$ and  $\mathcal K_i^1=\big\{ \{f_i^1;A^1\}\big\}$. Since each vertex of $A^2$ has colour 2 then each vertex of $F^1$ must have colour 1 or else a monochromatic $e$-star would now exist. Likewise, each vertex of $F^2$ must have colour 2. Let $\mathcal K^1=\bigcup\limits_{i=1}^{e}\mathcal K_i^1$ and $\mathcal K^2=\bigcup\limits_{i=1}^{e}\mathcal K_i^2$. We decompose the edges between $A^i$ and $F^i$, $1\leq i\leq 2$, later. Decompose the edges between $A^1$ and $G^2$ and between $A^2$ and $G^1$ into a set of $e$-stars $\mathcal L$ in a similar manner, forcing the vertices of the subsets $G^1$ and $G^2$ to have colours 1 and 2 respectively. We decompose the edges between $A^i$ and $G^i$, $1\leq i\leq 2$, later. Also, decompose the edges between $K^1$ and $G^2$ and between $K^2$ and $G^1$ into a set of $e$-stars $\mathcal C$ in a similar manner, forcing the vertices of the subsets $K^1$ and $K^2$ to have colours 1 and 2 respectively. We decompose the edges between $K^i$ and $G^i$, $1\leq i\leq 2$, later

Next, we decompose the edges between $H$ and $F\cup G$. For each $i\in \{1,\ldots,e\}$, let $\mathcal E_i^1=\big\{ \{h_i^1;f_1^2,\ldots,f_e^2\}$, $\{h_i^1;f_1^1,\ldots,f_{e-1}^1,g_1^2\}$, $\{h_i^1;f_e^1,g_1^1,\ldots,g_{e-2}^1,g_2^2\}$, $\{h_i^1;g_{e-1}^1,g_{e}^1,g_3^2,\ldots,g_e^2\}\big\}$, and $\mathcal E^1=\bigcup\limits_{i=1}^{e}\mathcal E_i^1$.  Then $\mathcal E^1$ is a decomposition of the edges between $H^1$ and $F\cup G$ into $e$-stars such that all the vertices in $H^1$ are forced to have colour 1. Decompose the edges between $H^2$ and $F\cup G$ into a set of $e$-stars $\mathcal E^2$ in a similar manner, forcing every vertex in $H^2$ to have unique colour $2$. Let $\mathcal E=\bigcup\limits_{j=1}^{2}\mathcal E^j$.

Partition the $6e-6$ vertices of $F^1\cup G^1\cup H\cup K\setminus \{h_1^2,h_2^2,h_3^2,k_1^2,k_2^2,k_3^2\}$ into  six sets $S_1,\ldots,S_6$ of size $e-1$. For each $1\leq i\leq e$, let $M_i^1=\big\{ \{a_i^1;S_1\cup \{h_1^2\}\},$ $\{a_i^1;S_2\cup \{h_2^2\}\},$ $\{a_i^1;S_3\cup \{h_3^2\}\},$ $\{a_i^1;S_4\cup \{k_1^2\}\},$ $\{a_i^1;S_5\cup \{k_2^2\}\},$ $\{a_i^1;S_6\cup \{k_3^2\}\}\big\}$. Then $M^1=\bigcup\limits_{i=1}^e M_i^1$ is a decomposition of the edges between $A^1$ and $F^1\cup G^1\cup H\cup K$ into $e$-stars, none of which is monochromatic. Decompose the edges between $A^2$ and $F^2\cup G^2\cup H\cup K$ into a set of $e$-stars $M^2$ in a similar manner. Let $\mathcal M=\bigcup\limits_{j=1}^{2}\mathcal M^j$.

Partition the $5e-5$ vertices of $G^1\cup F\cup H\setminus \{f_1^2,f_2^2,f_3^2,h_1^2,h_2^2\}$ into five sets $R_1,\ldots,R_5$ of size $e-1$. For each $1\leq i\leq e$, let $\mathcal F_i^1=\big\{ \{k_i^1;R_1\cup \{f_1^2\}\},$ $\{k_i^1;R_2\cup \{f_2^2\}\},$ $\{k_i^1;R_3\cup \{f_3^2\}\},$
 $\{k_i^1;R_4\cup \{h_1^2\}\},$ $\{k_i^1;R_5\cup \{h_2^2\}\}$. Then $\mathcal F^1=\bigcup\limits_{i=1}^e \mathcal F_i^1$ is a decomposition of the edges between $K^1$ and $G^1\cup F\cup H$ into $e$-stars, none of which is monochromatic. Decompose the edges between $K^2$ and $G^2\cup F\cup H$ into a set of $e$-stars $\mathcal F^2$ in a similar manner. 

Next, for each $1\leq i\leq \ell$, we decompose the edges between $F\cup G$ and $U_i$ into $e$-stars. Let $u_1^s,\ldots,u_{e}^s$ be the vertices of colour $s$, $1\leq s\leq 2$, in an equitable  $2$-colouring of $(U_0,\mathcal A_0)$. For each $j\in\{1,\ldots,e\}$, let $\mathcal N_j^1=\big\{\{u_j^1\times \{i\};f_1^2,\ldots,f_e^2\}$, $\{u_j^1\times \{i\};f_1^1,\ldots,f_{e-1}^1,g_1^2\}$, $\{u_j^1\times \{i\};f_e^1,g_1^1,\ldots,g_{e-2}^1,g_2^2\}$, $\{u_j^1\times \{i\}$;$g_{e-1}^1,g_e^1,$ $g_3^2,\ldots,g_e^2\}\big\}$ and observe that $u_j^1 \times \{i\}$ in now forced to have colour 1. 

For each $j\in\{1,\ldots,e\}$, decompose the edges between $u_j^2\times \{i\}$ and $F\cup G$ into a set of $e$-stars $\mathcal N_j^2$ in a similar manner so that $u_j^2\times \{i\}$ is forced to have colour 2. Let $\mathcal N^s=\bigcup\limits_{j=1}^{e} \mathcal N_j^s$ and $\mathcal O_i=\bigcup\limits_{s=1}^{2} \mathcal N^s$. Then $\mathcal O_i$ is a decomposition of the edges between $F\cup G$ and $U_i$ into $e$-stars which forces every vertex of $U_i$ to have a unique colour. Let $\mathcal O=\bigcup\limits_{i=1}^{\ell}\mathcal O_i$. 

At this point every vertex of $V$ is coloured, and the colouring is unique (up to a permutation of the colours). Moreover, the two colour classes are equal in size and so the colouring is strongly equitable. All that remains is to complete the decomposition without introducing monochromatic $e$-stars.

For each $1\leq i\leq e$, we begin to decompose the edges between $F$ and $G$ into $e$-stars. Let $\mathcal Q_i^1=\big\{\{f_i^1;g_1^1,\ldots,g_{e-1}^1,g_1^2\}$, $\{f_i^1;g_e^1,g_2^2,\ldots,g_e^2\}\big\}$. Then $\mathcal Q^1=\bigcup\limits_{i=1}^{e}\mathcal Q_i^1$ is a decomposition of the edges between $F^1$ and $G$ into $e$-stars, none of which is monochromatic. Decompose the edges between $F^2$ and $G$ into a set of $e$-stars $\mathcal Q^2$ in a similar manner. Let $\mathcal Q=\bigcup\limits_{j=1}^{2}\mathcal Q^j$.

Next, for each $1\leq i\leq \ell$, we decompose the edges between $H$ and $U_i$ into $e$-stars. Recall that $u_1^s,\ldots,u_e^s$ are vertices of colour $s$, $1\leq s\leq 2$. Let $\mathcal I_j^1=\big\{\{u_j^1\times \{i\};h_1^1,\ldots,h_{e-1}^1,h_1^2\}$, $\{u_j^1\times \{i\};h_e^1,h_2^2,\ldots,h_e^2\}\big\}$, $1\leq j\leq e$. For each $j\in\{1,\ldots,e\}$, decompose the edges between $u_j^2\times \{i\}$ and $H$ into a set of $e$-stars $\mathcal I_j^2$ in a similar manner. Let $\mathcal I^s=\bigcup\limits_{j=1}^{e}\mathcal I_j^s$ and $\mathcal S_i=\bigcup\limits_{s=1}^{2}\mathcal I^s$. Then $\mathcal S_i$ is a decomposition of the edges between $H$ and $U_i$ into $e$-stars. Let $\mathcal S=\bigcup\limits_{i=1}^{\ell}\mathcal S_i$. For each $1\leq i\leq \ell$, decompose the edges between $K$ and $U_i$ into a set of $e$-stars $\mathcal H_i$ in a similar manner so that $\bigcup\limits_{i=1}^{\ell} \mathcal H_i$ decomposes all of the edges between $K$ and $U$. The edges between $U_i$ and $U_j$, $1\leq i < j \leq \ell$, are decomposed into set of $e$-stars $\mathcal U_j^i$ in a similar manner.

Note that $|A|=|F|=|G|=|H|=|K|=2e$ and so we let each of $(A,\mathcal A),(F,\mathcal F),(G,\mathcal G)$, $(H,\mathcal H)$ and $(K,\mathcal K)$ be an equitable 2-chromatic $e$-star system on the sets $A$, $F$, $G$, $H$ and $K$ respectively such that the colour classes are in agreement with the colouring that we have forced upon $V$.

Finally, let $\hat{\mathcal B}=\mathcal A \cup \mathcal F \cup \mathcal G \cup \mathcal H \cup \mathcal K \cup \big( \bigcup\limits_{i=1}^{\ell} \mathcal A_i\big) \cup \big( \bigcup\limits_{i=1}^{\ell} \mathcal P_i\big) \cup \big( \bigcup\limits_{j=2}^{\ell} \mathcal U_j^1 \big) \cup \big( \bigcup\limits_{j=3}^{\ell} \mathcal U_j^2 \big)\cup \ldots \cup \big( \mathcal U_{\ell}^{\ell -1} \big) \cup \big(\bigcup\limits_{j=1}^{2} \mathcal K^j\big) \cup \mathcal L \cup \mathcal E\cup \mathcal M\cup \mathcal O\cup \mathcal Q\cup \mathcal S\cup \mathcal C\cup \big(\bigcup\limits_{j=1}^{2} \mathcal F^j \big)\cup \big(\bigcup\limits_{i=1}^{\ell} \mathcal H_i\big)$.
Then $(\hat{V},\hat{\mathcal B})$ is an $e$-star system of order $n=10e+\ell n_0$ which is strongly equitable uniquely $2$-chromatic.
\end{Proof}
Observe that each of the cardinalities of the colour classes of the uniquely $2$-chromatic $e$-star systems constructed from Theorem~\ref{thm4.1} is greater than $e$. We tacitly use this property to construct a uniquely $2$-chromatic $e$-star system from a smaller uniquely $2$-chromatic $e$-star system.
\begin{Theorem}\label{thm4.2}
Let $(V,\mathcal B)$ be a strongly equitable uniquely 2-chromatic $e$-star system of order $n_0$ constructed from Theorem~\ref{thm4.1} with colour classes $C_1$ and $C_2$. Then there exists a uniquely $2$-chromatic $e$-star system for all $n >n_0$ such that $n\equiv 0,1$ (mod $2e$).
\end{Theorem}
\begin{Proof}
For $s\in \{1,2\}$, let $C_s=\{c_1^s,\ldots,c_{|C_s|}^s\}$.
Since $n_0\equiv 0$ (mod $2e$), then let $n_0=2et$, $t\geq 1$.

First, we construct a uniquely $2$-chromatic $S_e(2et+1)$, $(\hat{V}, \hat{\mathcal B})$, from $(V,\mathcal B)$. Let $\hat{V}=V\cup \{2et+1\}$. Let $\mathcal T_1^{2et+1}=\big\{\{2et+1;c_1^2,\ldots,c_e^2\}\big\}$ and $\mathcal T_2^{2et+1}=\big\{\{2et+1;c_1^1,\ldots,c_{e-1}^1,c_{e+1}^2\}$, $\{2et+1;c_e^1,\ldots,c_{2e-2}^1,c_{e+2}^2\}$, $\ldots$, $\{2et+1;c_{(t-1)(e-1)+1}^1,\ldots,c_{t(e-1)}^1,c_{e+t}^2\}$, $\{2et+1;c_{t(e-1)+1}^1,\ldots,$ $c_{|C_1|}^1,$ $c_{e+t+1}^2,\ldots,c_{e+t+r}^2\}\big\}$ where $t=\lfloor \frac{|C_1|}{e-1}\rfloor$ and $r=e-(|C_1|-t(e-1))$. 

Partition the set $C_2\setminus \{c_1^2,\ldots,c_e^2,c_{e+1}^2,c_{e+2}^2,\ldots,c_{e+t+r}^2\}$ into a set $\mathbb E$ of $m$ disjoint $e$-subsets $X_1,\ldots,X_m$ where $m=\frac{|C_2|-(e+t+r)}{e}$. Let $\mathcal T_3^{2et+1}=\bigcup\limits_{i=1}^{m} \big\{\{2et+1;X_i\}\big\}$.

Let $\hat{\mathcal B} =\mathcal B\cup \mathcal T_1^{2et+1}\cup \mathcal T_2^{2et+1}\cup \mathcal T_3^{2et+1}$. Then $(\hat{V},\hat {\mathcal B})$ is a uniquely $2$-chromatic $e$-star system of order $2et+1$ with colour classes $C_1\cup \{2et+1\}$ and $C_2$.

 Next, we construct a uniquely $2$-chromatic $S_e(2et+2e)$, $(\mathring{V}, \mathring 
 {B})$, from $(\hat{V},\hat{\mathcal B})$. Let $\mathring{V}=\hat{V}\cup V'$ where $V'=\{2et+2,2et+3,\ldots,2et+2e\}$ and let $v_0 \in C_1 \setminus \{c_1^1,\ldots,c_e^1\}$. Let $(V'\cup v_0,\mathcal B')$ be a 2-chromatic $e$-star system of order $2e$ constructed from Theorem~\ref{thm4}. Let $\mathcal T_1^{i}=\bigcup\limits_{i}\big\{\{i;c_1^2,\ldots,c_e^2\}\big\}$, $i\in \{2et+3,2et+5,\ldots, 2et+2e-1\}$ and $ \mathcal T_1^{j}=\bigcup\limits_{j}\big\{\{j;c_1^1,\ldots,c_e^1\}\big\}$, $j\in \{2et+2,2et+4,\ldots,2et+2e\}$. For any $k \in \{2et+2,2et+3,\ldots,2et+2e\}$, decompose the remaining edges between vertex $k$ and set $V\setminus \{v_0\}$ into sets of $e$-stars $\mathcal T_2^k$ and $\mathcal T_3^k$ in a manner similar to the construction of $(\hat{V},\hat{\mathcal B})$ from $(V,\mathcal B)$.
 
 Let $\mathring{\mathcal B}=\hat{\mathcal B}\cup \mathcal B'\cup \big( \bigcup\limits_{k=2et+2}^{2et+2e}\mathcal T_1^k\big) \cup \big(\bigcup\limits_{k=2et+2}^{2et+2e} \mathcal T_2^k\big) \cup \big(\bigcup\limits_{k=2et+2}^{2et+2e} \mathcal T_3^k\big)$. Then $(\mathring{V},\mathring{\mathcal B})$ is a uniquely $2$-chromatic $e$-star system of order $2et+2e$ with colour classes $C_1\cup \{2et+3,2et+5,\ldots,2et+2e-1\}$ and $C_2 \cup \{2et+2,2et+4,\ldots,2et+2e\}$. 
\end{Proof}

We then obtain the following corollary.
\begin{Corollary}
	There exists some integer $n_0$ where $n_0 \equiv 0$ (mod $2e$) such that for all admissible $n \geq n_0$ where $n\equiv 0,1$ (mod $2e$), there exists a uniquely $2$-chromatic $e$-star system of order $n$ for any $e\geq 3$.
\end{Corollary}

\begin{Proof}
	Apply Theorem~\ref{thm4.1} and Theorem~\ref{thm4.2}.
\end{Proof}
We now show how to construct a strongly equitable $k$-chromatic $e$-star system from a strongly equitable uniquely $(k-1)$-chromatic $e$-star system.
\begin{Theorem}\label{thm4.3}
Let $k\geq 3$ and $e\geq 3$. If there exists a strongly equitable uniquely $(k-1)$-chromatic $e$-star system of order $n_{k-1}\equiv 0$ (mod $2e$) with colour classes $C_1,\ldots,C_{k-1}$ such that $|C_i|>e$, $1\leq i\leq k-1$, then there exists a strongly equitable $k$-chromatic $e$-star system of order $n_k \equiv 0$ (mod $2e$).
\end{Theorem}
\begin{Proof}
Let $(U_0,\mathcal A_0)$ be a strongly equitable uniquely $(k-1)$-chromatic $e$-star system of order $n_{k-1}$ with colour classes $C_1,\ldots,C_{k-1}$ such that $n_{k-1}\equiv 0$ (mod $2e$), $|C_1|= \cdots =|C_{k-1}|=r>e$ where $r$ is a positive integer and each vertex of $C_s$ has colour $s$, for $1\leq s\leq k-1$. For each $s\in \{1,\ldots,k-1\}$, let $C_s=\{c_1^s,\ldots,c_r^s\}$. For each $i\in \{1,\ldots,k\}$, let $U_i=U_0 \times \{i\}$ and $\mathcal A_i = \mathcal A_0\times \{i\}$ where $\mathcal A_0\times \{i\}$ denotes the set $\{S\times \{i\} \mid S\in \mathcal A_0\}$. So $(U_i,\mathcal A_i)$ has colour classes $C_1\times\{i\},\ldots,C_{k-1}\times\{i\}$. We construct a strongly equitable $k$-chromatic $e$-star system $(V,\mathcal B)$ where $V=\bigcup \limits_{i=1}^{k} U_i$. We need to decompose the edges between the subsets of $V$ into $e$-stars such that the resulting $e$-star system $(V,\mathcal B)$ is strongly equitably $k$-chromatic.

First, we decompose the edges between $U_1$ and $U_2$. For each $s\in \{1,\ldots,k-1\}$, let $D_s=\big\{c_1^s\times \{1\},\ldots,c_e^s\times \{1\}\big\}$. Let $\mathcal F_1^1=\bigcup\limits_{j=1}^{r} \big\{\{c_j^1\times \{2\};D_1\}, \{c_j^1\times \{2\};D_2\},\ldots,\{c_j^1\times \{2\};D_{k-1}\}\big\}$. Partition $U_1 \setminus \bigcup\limits_{s=1}^{k-1} D_s$ into $m=\frac {n_{k-1}-e(k-1)}{e}$ sets $G_1,\ldots,G_m$ such that for each $\ell \in\{1,\ldots,m\}$, the vertices of $G_{\ell}$ do not all have the same colour. Let $\mathcal F_1^2=\bigcup\limits_{j=1}^{r} \big\{\{c_j^1\times \{2\};G_1\}, \{c_j^1\times \{2\};G_2\},\ldots,\{c_j^1\times \{2\};G_m\}\big\}$ and $\mathcal F_1=\mathcal F_1^1\cup \mathcal F_1^2$. Then $F_1$ is a decomposition of the edges between $C_1\times \{2\}$ and $U_1$ into $e$-stars. In any putative $(k-1)$-colouring of the $e$-star system that we are building, the colouring of $U_1$ is unique since $(U_1,\mathcal A_1)$ is uniquely $(k-1)$-chromatic. The uniqueness of the $(k-1)$-colouring of $(U_1,\mathcal A_1)$ is such that all vertices of $D_s$ $(1\leq s\leq k-1)$ have colour $s$. But now some star of $\mathcal F_1^1$ is monochromatic if any vertex of $C_1\times \{2\}$ has colour $s$ for some $s\in \{1,\ldots,k-1\}$. Therefore the system that we are building is not $(k-1)$-colourable. By demonstrating a $k$-colouring, we establish that it is $k$-chromatic. Without loss of generality, actually colour each vertex of $D_s$ with colour $s$ and each vertex of $C_1\times \{2\}$ with colour $k$. 

Let $\mathcal F_2^1=\bigcup\limits_{j=1}^{r} \big\{\{c_j^2\times \{2\};D_1\}, \{c_j^2\times \{2\};D_3\}, \{c_j^2\times \{2\};D_4\}, \ldots, \{c_j^2\times \{2\};D_{k-1}\}\big\}$. Partition $U_1\setminus \big( D_1\cup (\bigcup\limits_{s=3}^{k-1}D_s)\big)$ into $m=\frac{n_{k-1}-e(k-2)}{e}$ sets $G_1,\ldots,G_m$ such that for each $\ell\in\{1,\ldots,m\}$, the vertices of $G_{\ell}$ do not all have the same colour. Let $\mathcal F_2^2=\bigcup\limits_{j=1}^{r}\big\{\{c_j^2\times \{2\};G_1\},\ldots,\{c_j^2\times \{2\};G_m\}\big\}$ and $\mathcal F_2=\mathcal F_2^1\cup \mathcal F_2^2$. Then $\mathcal F_2$ is a decomposition of the edges between $C_2\times \{2\}$ and $U_1$ into $e$-stars such that all the vertices in $C_2\times \{2\}$ are forced to have colour 2. For each $s\in \{3,\ldots,k-1\}$, decompose the edges between $C_s\times \{2\}$ and $U_1$ into a set of $e$-stars $\mathcal F_s$ in a similar manner such that all the vertices in $C_s\times \{2\}$ are forced to have colour $s$.

 Let $\mathcal U_2^1=\bigcup\limits_{s=1}^{k-1}\mathcal F_s$. Then $\mathcal U_2^1$ is a decomposition of the edges between $U_1$ and $U_2$ into $e$-stars forcing $(U_2,\mathcal A_2)$ to have colour classes $C_2^2,C_2^3,\ldots,C_2^k$ of equal size with no vertex of colour 1, where $C_2^2=C_2\times \{2\}$, $C_2^3=C_3\times \{2\}$, $\ldots$, $C_2^k=C_1\times \{2\}$. For each $i\in \{3,\ldots,k\}$, decompose the edges between $U_1$ and $U_i$ into a set of $e$-stars $\mathcal U_i^1$ in a similar manner forcing $(U_i,\mathcal A_i)$ to have colour classes $C_i^1, C_i^2,\ldots,C_i^{i-2},C_i^{i},\ldots,C_i^k$ of equal size with no vertex of colour $i-1$, where $C_i^1=C_1\times \{i\}$, $C_i^2=C_2\times \{i\}$, $\ldots$, $C_i^{i-2}=C_{i-2}\times \{i\}$, $C_i^i=C_i\times \{i\}$, $\ldots$, $C_i^k=C_{i-1}\times \{i\}$.

Next, for each $j\in\{3,\ldots,k\}$, we decompose the edges between $U_2$ and $U_j$. Partition $U_2$ into $m=\frac{n_{k-1}}{e}$ sets $G_1,\ldots,G_m$ such that for each $\ell \in \{1,\ldots,m\}$, the vertices of $G_{\ell}$ do not all have the same colour. For each $j\in \{3,\ldots,k\}$, let $\mathcal U_j^2=\bigcup\limits_{u\in U_j} \big\{\{u;G_1\},\{u;G_2\},\ldots,\{u;G_m\}\big\}$. Then $\mathcal U_j^2$ is a decomposition of the edges between $U_2$ and $U_j$ into $e$-stars such that no $e$-star is monochromatic. For each $3\leq i<j\leq k$, decompose the edges between $U_i$ and $U_j$ into a set of $e$-stars $\mathcal U_j^i$ in a similar manner. 

Let $\mathcal B=\big(\bigcup\limits_{i=1}^{k}\mathcal A_i\big) \cup \big(\bigcup\limits_{j=2}^{k}\mathcal U_j^1\big) \cup \big(\bigcup\limits_{j=3}^{k}\mathcal U_j^2\big)\cup \ldots \cup \mathcal U_{k}^{k-1}$ and $C_1^1=C_1\times \{1\}$, $C_1^2=C_2\times \{1\}$, $\ldots$, $C_1^{k-1}=C_{k-1}\times \{1\}$. Then $(V,\mathcal B)$ is a strongly equitable $k$-chromatic $e$-star system of order $n_k=kn_{k-1}$ with colour classes $C_1^1\cup C_3^1\cup C_4^1\cup \cdots \cup C_k^1$,  $C_1^2\cup C_2^2\cup C_4^2\cup \cdots \cup C_k^2$, $\ldots$, $C_1^{k-1}\cup C_2^{k-1}\cup \cdots \cup C_{k-1}^{k-1}$, $C_2^k\cup C_3^k\cup \cdots \cup C_k^k$ where each colour class is of size $(k-1)r$.
\end{Proof}
We now show how to construct a strongly equitable uniquely $k$-chromatic $e$-star system from a strongly equitable $k$-chromatic $e$-star system.
\begin{Theorem}\label{thm4.4}
	Let $k\geq 3$ and $e\geq 3$. If there exists a strongly equitable $k$-chromatic $e$-star system of order $n_0\equiv 0$ (mod $2e$), then there exists a strongly equitable uniquely $k$-chromatic $e$-star system of order $n$ for some $n \equiv 0$ (mod $2e$).
\end{Theorem}
\begin{Proof}
	Let $(U_0,\mathcal A_0)$ be a strongly equitable $k$-chromatic $e$-star system of order $n_0\equiv 0$ (mod $2e$) with colour classes $C_1,\ldots,C_k$ such that each vertex of $C_s$ has colour $s$, for $s=1,\ldots,k$. For a positive integer $\ell$ that will be fixed later and for each $i\in \{1,\ldots,\ell\}$, let $U_i=U_0\times \{i\}$ and $\mathcal A_i=\mathcal A_0\times \{i\}$ where $\mathcal A_0 \times \{i\}$ denotes the set $\big\{ S\times \{i\} \, | \, S\in \mathcal A_0\big\}$. So $(U_i,\mathcal A_i)$ has colour classes $C_1\times \{i\},\ldots,C_k\times \{i\}$. Let $U=\bigcup\limits_{i=1}^{\ell} U_i$.
	Let $A^1=\{a_1^1,\ldots,a_{e}^1\},\ldots,A^k=\{a_1^k,\ldots,a_{e}^2k\}$, $A_0^1=\{(a_1)_0^1,\ldots,(a_e)_0^1\}$, $A=\bigcup_{i=1}^{k}A^i$ and $A'=A\cup A_0^1$. Similarly let $F^1=\{f_1^1,\ldots,f_{e}^1\},\ldots,F^k=\{f_1^k,\ldots,f_{e}^k\}$, $F_0^1=\{(f_1)_0^1,\ldots,(f_e)_0^1\}$, $F=\bigcup_{i=1}^{k}F^i$ and $F'=F\cup F_0^1$, $G^1=\{g_1^1,\ldots,g_{e}^1\},\ldots,G^k=\{g_1^k,\ldots,g_{e}^k\}$, $G_0^1=\{(g_1)_0^1,\ldots,(g_e)_0^1\}$, $G=\bigcup_{i=1}^{k}G^i$ and $G'=G\cup G_0^1,$ $H^1=\{h_1^1,\ldots,h_{e}^1\},\ldots,H^k=\{h_1^k,\ldots,h_{e}^k\}$, $H_0^1=\{(h_1)_0^1,\ldots,(h_e)_0^1\}$, $H=\bigcup_{i=1}^{k}H^i$ and $H'=H\cup H_0^1$. We construct a uniquely $k$-chromatic $e$-star system $(\hat{V},\hat{\mathcal B})$ where $\hat{V}=A\cup F\cup G\cup H\cup U$ if $k$ is even and $\hat{V}=A'\cup F'\cup G'\cup H'\cup U$ if $k$ is odd. 
	
	We need to decompose the edges between the subsets of $\hat{V}$ into $e$-stars such that the resulting $e$-star system $(\hat{V},\hat{\mathcal B})$ is uniquely $k$-chromatic. 
	
	First, we decompose the edges between $A$ and $U$ into $e$-stars in a way such that no $e$-subset in $A$ other than the $e$-subsets $A^1,\ldots,A^k$ are monochromatic in any putative $k$-colouring of $(\hat{V},\hat{\mathcal B})$. Let $\mathcal D$ be the set of all $e$-subsets of $A$ except 
	$A^1,\ldots,A^k$. The number of $e$-subsets of the set $A$ is $N=\binom{ke}{e}$ and $|\mathcal D|=N-k$. Let $a=k-1$ and  $|\mathcal D|=aq+r$ where $q$ is a nonnegative integer and $0\leq r< a$ is an integer. Partition set $\mathcal D$ into $q$ sets of disjoint $e$-subsets $\mathbb T_1,\ldots,\mathbb T_q$ of size $a$ and one set of disjoint $e$-subsets $\mathbb T_{q+1}$ of size $r$; such a partition is known to exist by Theorem~1 in \cite{Baranyai1975}. Let $\mathbb T_i=\{T_1^{i},\ldots,T_a^{i}\}$, $1\leq i\leq q$ and $\mathbb T_{q+1}=\{T_1^{q+1},\ldots,T_r^{q+1}\}$. Let $\tilde{A_i}=A\setminus \bigcup_{j=1}^{a} T_j^{i}$, $1\leq i\leq q$ and $\tilde{A}_{q+1}=A\setminus \bigcup_{j=1}^{r} T_j^{q+1}$ if $r>0$ and 
	$\tilde{A}_{q+1}=\emptyset$ if $r=0$. We now fix $\ell=q+1$ if $r>0$ and $\ell=q$ if $r=0$.
	
	Partition $U_0$ into $e$-subsets $E_1,\ldots,E_m$  where $m=\frac{n_0}{e}$ such that each $E_j$ has nonempty intersection with at least two colour classes, $1\leq j\leq m$. For each $i\in \{1,2,\ldots,\ell\}$ and for each $w\in \tilde {A_i}$, we decompose the edges between $w$ and $U_i$ into the set of $e$-stars $\mathcal R_w=\bigcup\limits_{j=1}^{m}\big\{ \{w;E_j\times \{i\} \}\big\}$. For each $i\in \{1,\ldots,q\}$, let $\mathcal T_i=\bigcup \limits_{u\in U_i} \big\{ \{u;T_1^i\},\ldots,\{u;T_a^{i}\big\}$ and $\mathcal T_{q+1}=\bigcup \limits_{u\in U_i} \big\{ \{u;T_1^{q+1}\},\ldots,\{u;T_r^{q+1}\}\big\}$. Then $\mathcal P_i=\mathcal T_i\cup \big (\bigcup \limits_{w\in \tilde {A_i}} \mathcal R_w \big)$, $1\leq i\leq \ell$,  is a decomposition of the edges between $A$ and $U_i$ and $\bigcup\limits_{i=1}^{\ell} \mathcal P_i$ is a decomposition of the edges between $A$ and $U$ into $e$-stars.
	
	Recall that $(U_i,\mathcal A_i)$ is $k$-chromatic. So for any $k$-colouring of $(U_i,\mathcal A_i)$ there must be a vertex $u_s^i$ of each colour $s$ in $U_i$, $1\leq s\leq k$. Note that for each $i\in \{1,\ldots,q\}$, $\{u_s^i;T_j^i\}\in \mathcal T_i$, $1\leq j\leq a$, and $\{u_s^{q+1};T_j^{q+1}\}\in \mathcal T_{q+1}$, $1\leq j\leq r$, and so to ultimately obtain a valid $k$-colouring of $(\hat{V},\hat{\mathcal B})$, the $e$ vertices of $T_j^i$ and $T_j^{q+1}$ must not all have colour $s$, $1\leq s\leq k$. If a colour $s$ occurs on more than $e$ vertices of $A$ then there exists a $T_j^i$ or $T_j^{q+1}$ with only colour $s$ and then $\{u_s^i;T_j^i\}$ or $\{u_s^i;T_j^{q+1}\}$ would be monochromatic. So in any valid $k$-colouring of $(\hat{V},\hat{\mathcal B})$, each colour must occur on exactly $e$ points of $A$. If the points of colour $s$ are not all in $A^1$ or $\ldots$ or $A^k$ then they are the leaves of a monochromatic $e$-star with centre $u_s^i$. Without loss of generality we therefore assign colour $s$ to all $e$ points of $A^s$, $1\leq s\leq k$.

	We now begin to decompose the edges between $A$ and $F$ into $e$-stars. For each $i\in \{1,\ldots,e\}$, let $\mathcal K_i^1=\big\{ \{f_i^1;A^2\},\ldots,\{f_i^1;A^k\}\big\}$, $\mathcal K_i^2=\big\{ \{f_i^2;A^1\},\{f_i^2;A^3\},\ldots,\{f_i^2;A^{k}\}\big\}$, $\ldots$, $\mathcal K_i^k=\big\{ \{f_i^k;A^1\},\ldots,\{f_i^k;A^{k-1}\}\big\}$. Since each vertex of $A^s$ has colour $s$, $2\leq s\leq k$, then each vertex of $F^1$ must have colour 1 or else a monochromatic $e$-star would now exist. Likewise, each vertex of $F^s$ must have colour $s$, $2\leq s\leq k$. Let $\mathcal K^1=\bigcup\limits_{i=1}^{e}\mathcal K_i^1,\ldots,\mathcal K^k=\bigcup\limits_{i=1}^{e}\mathcal K_i^k$. We decompose the edges between $A^i$ and $F^i$, $1\leq i\leq k$, later. For each $j\in \{1,\ldots,k\}$, decompose the edges between $A^j$ and $G^1\cup\cdots\cup G^{j-1}\cup G^{j+1}\cup \cdots\cup G^{k}$ into a set of $e$-stars $\mathcal L_j$ in a similar manner, forcing the vertices of the subsets $G^1,\ldots,G^k$ to have colours $1,\ldots,k$ respectively. Let $\mathcal L=\bigcup\limits_{j=1}^{k} \mathcal L_j$. We decompose the edges between $A^i$ and $G^i$, $1\leq i\leq k$, later.
	
	Next, we decompose the edges between $H$ and $F\cup G$. For each $i\in \{1,\ldots,e\}$, let $\mathcal E_i^1=\big\{ \{h_i^1;f_1^2,\ldots,f_e^2\},\ldots,$ $\{h_i^1;f_1^k,\ldots,f_e^k\}$, $\{h_i^1;f_1^1,\ldots,f_{e-1}^1,g_1^2\}$, $\{h_i^1;f_e^1,g_1^1,\ldots,g_{e-2}^1,g_2^2\}$, $\{h_i^1;g_{e-1}^1,$ $g_{e}^1,g_3^2,\ldots,g_e^2\}$, $\{h_i^1;g_1^3,\ldots,g_{e}^3\},\ldots,$ $\{h_i^1;g_1^k,\ldots,g_{e}^k\}\big\}$ and $\mathcal E^1=\bigcup\limits_{i=1}^{e}\mathcal E_i^1$.  Then $\mathcal E^1$ is a decomposition of the edges between $H^1$ and $F\cup G$ into $e$-stars such that all the vertices in $H^1$ are forced to have colour 1. For each $s\in \{2,\ldots,k\}$, decompose the edges between $H^s$ and $F\cup G$ into a set of $e$-stars $\mathcal E^s$ in a similar manner, forcing every vertex in $H^s$ to have unique colour $s$. Let $\mathcal E=\bigcup\limits_{s=1}^{k}\mathcal E^s$.

	Next, we decompose the edges between $A^1$ and $F^1\cup G^1\cup H$. For each $i\in \{1,\ldots,e\}$, let $\mathcal M_i^1=\big\{\{a_i^1;f_1^1,\ldots,f_{e-1}^1,h_1^2\},$ $\{a_i^1;f_e^1,g_1^1,\ldots,g_{e-2}^1,h_2^2\},$ $\{a_i^1;g_{e-1}^1,g_e^1,h_3^2,\ldots,h_e^2\},$
	$\{a_i^1;h_1^1,\ldots,h_{e-1}^1,h_1^3\},$
	$\{a_i^1;h_e^1,h_2^3,\ldots,h_e^3\},$
	$\{a_i^1;h_1^4,\ldots,h_e^4\},\ldots,$
	$\{a_i^1;h_1^k,\ldots,h_e^k\}\big\}$. Then $\mathcal M^1=\bigcup\limits_{i=1}^{e} \mathcal M_i^1$ is a decomposition of the edges between $A^1$ and $F^1\cup G^1\cup H$ into $e$-stars, none of which is monochromatic. For each $j\in \{2,\ldots,k\}$, decompose the edges between $A^j$ and $F^j\cup G^j\cup H$ into a set of $e$-stars $\mathcal M^j$ in a similar manner. Let $\mathcal M=\bigcup\limits_{j=1}^{k}\mathcal M^j$.

	Next, for each $1\leq i\leq \ell$, we decompose the edges between $F\cup G$ and $U_i$ into $e$-stars. Let $u_1^s,\ldots,u_{e}^s$ be the vertices of colour $s$, $1\leq s\leq k$, in an equitable  $k$-colouring of $(U_0,\mathcal A_0)$. For each $j\in \{1,\ldots,e\}$, let $\mathcal N_j^1=\big\{\{u_j^1\times \{i\};f_1^2,\ldots,f_e^2\},\ldots,$ $\{u_j^1\times \{i\};f_1^k,\ldots,f_e^k\}$, $\{u_j^1\times \{i\};f_1^1,\ldots,f_{e-1}^1,g_1^2\}$, $\{u_j^1\times \{i\};f_e^1,g_1^1,\ldots,g_{e-2}^1,g_2^2\}$, $\{u_j^1\times \{i\}$;$g_{e-1}^1,g_e^1,$ $g_3^2,\ldots,g_e^2\}$, $\{u_j^1\times \{i\};g_1^3,$ $\ldots,g_e^3\},\ldots,$ $\{u_j^1\times \{i\};g_1^k,\ldots,g_e^k\}\big\}$ and observe that $u_j^1 \times \{i\}$ is now forced to have colour 1. For each $j\in\{1,\ldots,e\}$ and each $s\in \{2,\ldots,k\}$, decompose the edges between $u_j^s\times \{i\}$ and $F\cup G$ into a set of $e$-stars $\mathcal N_j^s$ in a similar manner so that $u_j^s\times \{i\}$ is forced to have colour $s$. Let $\mathcal N^s=\bigcup\limits_{j=1}^{e} \mathcal N_j^s$ and $\mathcal O_i=\bigcup\limits_{s=1}^{k} \mathcal N^s$. Then $\mathcal O_i$ is a decomposition of the edges between $F\cup G$ and $U_i$ into $e$-stars which forces every vertex of $U_i$ to have a unique colour. Let $\mathcal O=\bigcup\limits_{i=1}^{\ell}\mathcal O_i$. 
	
	At this point every vertex of $\hat{V}$ is coloured, and the colouring is unique (up to a permutation of the colours). Moreover, the colour classes are equal in size and so the colouring is strongly equitable. All that remains is to complete the decomposition without introducing monochromatic $e$-stars.
	
	For each $1\leq i\leq e$, we begin to decompose the edges between $F$ and $G$ into $e$-stars. Let $\mathcal Q_i^1=\big\{\{f_i^1;g_1^1,\ldots,g_{e-1}^1,g_1^2\}$, $\{f_i^1;g_e^1,g_2^2,\ldots,g_e^2\}\big\}$, $\{f_i^1;g_1^3,\ldots,g_e^3\},\ldots,$ $\{f_i^1;g_1^k,\ldots,g_e^k\}\big\}$. Then $\mathcal Q^1=\bigcup\limits_{i=1}^{e}\mathcal Q_i^1$ is a decomposition of the edges between $F^1$ and $G$ into $e$-stars, none of which is monochromatic. For each $j\in \{2,\ldots,k\}$, decompose the edges between $F^j$ and $G$ into a set of $e$-stars $\mathcal Q^j$ in a similar manner. Let $\mathcal Q=\bigcup\limits_{j=1}^{k}\mathcal Q^j$.
	
	Next, for each $1\leq i\leq \ell$, we decompose the edges between $H$ and $U_i$ into $e$-stars. Recall that $u_1^s,\ldots,u_e^s$ are the vertices of colour $s$, $1\leq s\leq k$. Let $\mathcal I_j^1=\big\{\{u_j^1\times \{i\};h_1^1,\ldots,h_{e-1}^1,h_1^2\}$, $\{u_j^1\times \{i\};h_e^1,h_2^2,\ldots,h_e^2\}$, $\{u_j^1\times \{i\};h_1^3,\ldots,h_e^3\},\ldots,$ $\{u_j^1\times\{i\};h_1^k,\ldots,h_e^k\}\big\}$, $1\leq j\leq e$. For each $j\in\{1,\ldots,e\}$ and each $s\in \{2,\ldots,k\}$, decompose the edges between $u_j^s\times \{i\}$ and $H$ into a set of $e$-stars $\mathcal I_j^s$ in a similar manner. Let $\mathcal I^s=\bigcup\limits_{j=1}^{e}\mathcal I_j^s$ and $\mathcal S_i=\bigcup\limits_{s=1}^{k}\mathcal I^s$. Then $\mathcal S_i$ is a decomposition of the edges between $H$ and $U_i$ into $e$-stars. Let $\mathcal S=\bigcup\limits_{i=1}^{\ell}\mathcal S_i$. The edges between $U_i$ and $U_j$, $1\leq i < j \leq \ell$, are decomposed into set of $e$-stars $\mathcal U_j^i$ in a similar manner.
	
	Note that $|A|=|F|=|G|=|H|=ke$.
	
	\textbf{Case1.} $k$ is even. We let each of $(A,\mathcal A),(F,\mathcal F),(G,\mathcal G)$, and $(H,\mathcal H)$ be strongly equitable $k$-chromatic $e$-star systems on the sets $A$, $F$, $G$, and $H$ respectively such that the colour classes are in agreement with the colouring that we have forced upon $\hat{V}$. Let $\hat{\mathcal B}=\mathcal A \cup \mathcal F \cup \mathcal G \cup \mathcal H \cup \big( \bigcup\limits_{i=1}^{\ell} \mathcal A_i\big) \cup \big( \bigcup\limits_{i=1}^{\ell} \mathcal P_i\big) \cup \big( \bigcup\limits_{j=2}^{\ell} \mathcal U_j^1 \big) \cup \big( \bigcup\limits_{j=3}^{\ell} \mathcal U_j^2 \big)\cup \ldots \cup \big( \mathcal U_{\ell}^{\ell -1} \big) \cup \big(\bigcup\limits_{j=1}^{k} \mathcal K^j\big) \cup \mathcal L \cup \mathcal E\cup \mathcal M\cup \mathcal O\cup \mathcal Q\cup \mathcal S$. Then $(\hat{V},\hat{\mathcal B})$ is an $e$-star system of order $n=4ke+\ell n_0$ which is uniquely $k$-chromatic and strongly equitable.

\textbf{Case 2.} $k$ is odd. For $j \in \{1,\ldots,e\}$, let $\mathcal V_j=\big\{\{(a_j)_0^1;f_1^2,\ldots,f_e^2\},$
$\{(a_j)_0^1;f_1^3,\ldots,f_e^3\},$ $\ldots,$ $\{(a_j)_0^1;f_1^k,\ldots,f_e^k\},$ $\{(a_j)_0^1;f_1^1,\ldots,f_{e-1}^1,g_1^2\},$ $\{(a_j)_0^1;f_e^1,g_1^1,\ldots,g_{e-2}^1,g_2^2\},$ $\{(a_j)_0^1;g_{e-1}^1,g_e^1,g_3^2,\ldots,g_e^2\},$ $\{(a_j)_0^1;$ $g_1^3,$ $\ldots,g_e^3\},$ \ldots,$ \{(a_j)_0^1;g_1^k,\ldots,g_e^k\}\big\}$ $\cup$ $\big\{\{(a_j)_0^1;h_1^1,\ldots,h_{e-1}^1,h_1^2\},$ $\{(a_j)_0^1;h_e^1,h_2^2,\ldots,h_{e-1}^2\},$ $\{(a_j)_0^1;h_1^3,\ldots,h_e^3\},$ $\ldots,$ $\{(a_j)_0^1;h_1^k,\ldots,h_e^k\}\big\}$ $\cup$ $\bigcup\limits_{i=1}^{\ell}\big(\big\{\{(a_j)_0^1;u_1^1\times \{i\},u_2^1\times \{i\},\ldots,u_{e-1}^1\times \{i\},u_1^2\times \{i\}\},$ $\ldots,$ $\{(a_j)_0^1;u_{(t-1)(e-1)+1}^1\times \{i\},\ldots,u_{t(e-1)}^1\times \{i\},u_{t}^2\times \{i\}\}$, $\{(a_j)_0^1;u_{t(e-1)+1}^1\times \{i\},\ldots,u_{|C_1|}^1\times \{i\},u_{t+1}^2\times \{i\},\ldots,u_{t+m}^2\times \{i\}\}\big\}\big)$ $\cup$ $\bigcup\limits_{i=1}^{\ell}\big(\big\{\{(a_j)_0^1;x_1^d\times \{i\},\ldots,x_e^d\times \{i\}\}\mid 1\leq d\leq \frac{n_0-(|C_1|+t+m)}{e}\big\}\big)$, where $t$ and $m$ are positive integers and the sets $\{x_1^d,\ldots,x_e^d\}$, $1\leq d\leq \frac{n_0-(|C_1|+t+m)}{e}$, form a partition of $U_i\setminus (C_1\times \{i\}) \cup \big\{u_1^2\times \{i\},u_2^2\times \{i\},\ldots, u_{t+m}^2\times \{i\}\big\}$. Let $\mathcal V=\bigcup\limits_{j=1}^{e}\mathcal V_j$, then $\mathcal V$ is a decomposition of the edges between $A_0^1$ and $F\cup G\cup H\cup U$ into $e$-stars forcing the vertices of $A_0^1$ to have unique colour~1.

Decompose the edges between $F_0^1$ and $A'\cup G\cup H\cup U$ into set of $e$-stars $\mathcal W$ in a similar manner so that the vertices of $F_0^1$ are forced to have colour~1. Also, decompose the edges between $G_0^1$ and $A'\cup F'\cup H\cup U$ into set of $e$-stars $\mathcal X$ so that the vertices of $G_0^1$ are forced to have colour~1 and the edges between $H_0^1$ and $A'\cup F'\cup G'\cup U$ into set of $e$-stars $\mathcal Y$ so that the vertices of $H_0^1$ are forced to have colour~1.

Now, we decompose the edges of the complete graph on set $A'$ into $e$-stars. We first take 2-chromatic $e$-star systems $Z_1=(A^1\cup A^2,\mathcal A_2^1)$, $Z_2=(A^3\cup A^4,\mathcal A_4^3)$, $\ldots,$ $Z_{\frac{k+1}{2}}=(A^k\cup A_0^1,\mathcal A_1^k)$ of size $2e$. It is easy then to decompose the edges between $Z_i$ and $Z_j$ into set of $e$-stars $\mathcal Z_j^i$, $1\leq i<j\leq \frac{k+1}{2}$ such that no $e$-star is monochromatic. Then $(A',\mathcal A')$ where $\mathcal A'=\mathcal A_2^1\cup$ $\mathcal A_4^3\cup$ $\cdots$ $\cup$ $\mathcal A_1^k$ $\cup$ $(\bigcup\limits_{j=2}^{\frac{k+1}{2}}\mathcal Z_j^1)\cup$ $(\bigcup\limits_{j=3}^{\frac{k+1}{2}}\mathcal Z_j^2)\cup$ $\cdots \cup \mathcal Z_{\frac{k+1}{2}}^{\frac{k+1}{2}-1}$ is an $e$-star system which is uniquely $k$-chromatic. Decompose the edges of the complete graph on the sets $F'$, $G'$ and $H'$ into set of $e$stars $\mathcal F'$, $\mathcal G'$ and $\mathcal H'$ in a similar manner.

Let $\hat{\mathcal B}=\mathcal A' \cup \mathcal F' \cup \mathcal G' \cup \mathcal H' \cup \big( \bigcup\limits_{i=1}^{\ell} \mathcal A_i\big) \cup (\bigcup\limits_{i=1}^{\ell}\mathcal P_i) \cup \big( \bigcup\limits_{j=2}^{\ell} \mathcal U_j^1 \big) \cup \big( \bigcup\limits_{j=3}^{\ell} \mathcal U_j^2 \big)\cup \ldots \cup \big( \mathcal U_{\ell}^{\ell -1} \big) \cup (\bigcup\limits_{j=1}^{k}\mathcal K^j) \cup \mathcal L \cup \mathcal E\cup \mathcal M\cup \mathcal O\cup \mathcal Q\cup \mathcal S \cup \mathcal V\cup \mathcal W\cup \mathcal X\cup \mathcal Y$.
Then $(\hat{V},\hat{\mathcal B})$ is an $e$-star system of order $n=4e(k+1)+\ell n_0$ which is uniquely $k$-chromatic and strongly equitable.
\end{Proof}
We now show how to construct a uniquely $k$-chromatic $e$-star system from a smaller uniquely $k$-chromatic $e$-star system.
\begin{Theorem}\label{thm4.5}
	For any $k\geq 3$, let $(V,\mathcal B)$ be a strongly equitable uniquely $k$-chromatic $e$-star system of order $n_0$ constructed from Theorem~\ref{thm4.4} with colour classes $C_1,\ldots,C_k$. Then there exists a uniquely $k$-chromatic $e$-star system for all $n >n_0$ such that $n\equiv 0,1$ (mod $2e$).
\end{Theorem}
\begin{Proof}
	For $s\in \{1,\ldots,k\}$, let $C_s=\{c_1^s,\ldots,c_{|C_s|}^s\}$.
	Since $n_0\equiv 0$ (mod $2e$), then let $n_0=2et$, $t\geq 1$.
	
	First, we construct a uniquely $k$-chromatic $S_e(2et+1)$, $(\hat{V}, \hat{\mathcal B})$, from $(V,\mathcal B)$. Let $\hat{V}=V\cup \{2et+1\}$. Let $\mathcal T_1^{2et+1}=\big\{\{2et+1;c_1^2,\ldots,c_e^2\}, \{2et+1;c_1^3,\ldots,c_e^3\}, \ldots, \{2et+1;c_1^k,\ldots,c_e^k\} \big\}$ and $\mathcal T_2^{2et+1}=\big\{\{2et+1;c_1^1,\ldots,c_{e-1}^1,c_{e+1}^2\}$, $\{2et+1;c_e^1,\ldots,c_{2e-2}^1,c_{e+2}^2\}$, $\ldots$, $\{2et+1;c_{(t-1)(e-1)+1}^1,\ldots,$ $c_{t(e-1)}^1,c_{e+t}^2\}$, $\{2et+1;c_{t(e-1)+1}^1,\ldots,$ $c_{|C_1|}^1,$ $c_{e+t+1}^2,\ldots,c_{e+t+r}^2\}\big\}$ where $t=\lfloor \frac{|C_1|}{e-1}\rfloor$ and $r=e-(|C_1|-t(e-1))$. 
	
	Partition the set $(C_2\setminus \{c_1^2,\ldots,c_e^2,c_{e+1}^2,c_{e+2}^2,\ldots,c_{e+t+r}^2\})\cup (C_3\setminus \{c_1^3,\ldots,c_e^3\})\cup \cdots \cup (C_k\setminus \{c_1^k,\ldots,c_e^k\})$ into a set $\mathbb E$ of $m$ disjoint $e$-subsets $X_1,\ldots,X_m$ where
 $m=\frac{(|C_2|-(e+t+r))+(|C_3|-e)+\cdots+(|C_k|-e)}{e}$. Let $\mathcal T_3^{2et+1}=\bigcup\limits_{i=1}^{m} \big\{\{2et+1;X_i\}\big\}$.
	
	Let $\hat{\mathcal B} =\mathcal B\cup \mathcal T_1^{2et+1}\cup \mathcal T_2^{2et+1}\cup \mathcal T_3^{2et+1}$. Then $(\hat{V},\hat {\mathcal B})$ is a uniquely $k$-chromatic $e$-star system of order $2et+1$ with colour classes $C_1\cup \{2et+1\},C_2,\ldots,C_k$.

	Next, we construct a uniquely $k$-chromatic $S_e(2et+2e)$, $(\mathring{V}, \mathring 
	{B})$, from $(\hat{V},\hat{\mathcal B})$. Let $\mathring{V}=\hat{V}\cup V'$ where $V'=\{2et+2,2et+3,\ldots,2et+2e\}$ and let $v_0 \in C_1 \setminus \{c_1^1,\ldots,c_e^1\}$. Let $(V'\cup \{v_0\},\mathcal B')$ be a strongly equitable 2-chromatic $e$-star system of order $2e$ constructed from Theorem~\ref{thm4}. Let $\mathcal T_1^{i}=\bigcup\limits_{i}\big\{\{i;c_1^2,\ldots,c_e^2\},\{i;c_1^3,\ldots,c_e^3\},\ldots,\{i;c_1^k,\ldots,c_e^k\}\big\}$, $i\in \{2et+3,2et+5,\ldots, 2et+2e-1\}$ and $ \mathcal T_1^{j}=\bigcup\limits_{j}\big\{\{j;c_1^1,\ldots,c_e^1\},\{j;c_1^3,\ldots,c_e^3\},\ldots,\{j;c_1^k,\ldots,c_e^k\}\big\}$, $j\in \{2et+2,2et+4,\ldots,2et+2e\}$. For any $k \in \{2et+2,2et+3,\ldots,2et+2e\}$, decompose the remaining edges between vertex $k$ and set $V\setminus \{v_0\}$ into sets of $e$-stars $\mathcal T_2^k$ and $\mathcal T_3^k$ in a manner similar to the construction of $\mathcal T_2^{2et+1}$ and $\mathcal T_3^{2et+1}$ in $(\hat{V},\hat{\mathcal B})$.
	
	Let $\mathring{\mathcal B}=\hat{\mathcal B}\cup \mathcal B'\cup \big( \bigcup\limits_{k=2et+2}^{2et+2e}\mathcal T_1^k\big) \cup \big(\bigcup\limits_{k=2et+2}^{2et+2e} \mathcal T_2^k\big) \cup \big(\bigcup\limits_{k=2et+2}^{2et+2e} \mathcal T_3^k\big)$. Then $(\mathring{V},\mathring{\mathcal B})$ is a uniquely $k$-chromatic $e$-star system of order $2et+2e$ with colour classes $C_1\cup \{2et+3,2et+5,\ldots,2et+2e-1\},C_2 \cup \{2et+2,2et+4,\ldots,2et+2e\},C_3,\ldots,C_k$. 
\end{Proof}
We then conclude with the following corollary.
\begin{Corollary}
	There exists some integer $n_0$ where $n_0 \equiv 0$ (mod $2e$) such that for all admissible $n \geq n_0$ where $n\equiv 0,1$ (mod $2e$), there exists a uniquely $k$-chromatic $e$-star system of order $n$ for any $e\geq 3$ and $k\geq 3$.
\end{Corollary}

\begin{Proof}
	Apply Theorem~\ref{thm4.4} and Theorem~\ref{thm4.5}.
\end{Proof}

\section{Acknowledgements}
D.A.~Pike acknowledges research support from NSERC (grant number RGPIN-04456-2016).

\end{doublespace}

\end{document}